\magnification=\magstep1   
\input amstex
\UseAMSsymbols
\input pictex 
\vsize=23truecm
\NoBlackBoxes
\parindent=18pt
  
   \font\rmk=cmr8

\def\op{\text{\rm op}}

\def\mod{\operatorname{mod}}

\def\Hom{\operatorname{Hom}}
\def\End{\operatorname{End}}
\def\Ext{\operatorname{Ext}}

\def\Ann{\operatorname{Ann}}

\def\Ker{\operatorname{Ker}}

\def\soc{\operatorname{soc}}

\def\Im{\operatorname{Im}}

\def\bdim{\operatorname{\bold{dim}}}

\def\ss{\ssize}

\def\s{\hfill \square} 
\vglue1cm
\centerline{\bf The short local algebras of dimension 6}
	\smallskip
\centerline{\bf with non-projective reflexive modules.}
	\bigskip 
\centerline{Claus Michael Ringel}
	\bigskip\bigskip
\noindent {\narrower Abstract:  \rmk 
Let $\ss A$ be a finite-dimensional local 
algebra over an algebraically closed field, let $\ss J$ be the
radical of $\ss A.$ The modules we are interested in are the finitely generated left
$\ss A$-modules. Projective modules are always reflexive, and 
an algebra is self-injective iff all modules are reflexive.
We discuss the existence of non-projective reflexive modules in case $\ss A$ is
not self-injective. We assume that $\ss A$ is short (this means that $\ss J^3 = 0$).
In a joint paper with Zhang Pu, it has been shown that 6 is the smallest possible dimension 
of $\ss A$ that can occur and that in this case the following conditions have to be satisfied: 
$\ss J^2$ is both the left socle and the right 
socle of $\ss A$ and there is no uniform ideal of length 3. The present paper is devoted to show
the converse.

	\medskip
\noindent 
\rmk Key words. Short local algebra. Reflexive module. Gorenstein-projective module. 
Bristle. Atom. Bar. Bristle-bar layout.

\par}
	\bigskip\bigskip 

{\bf 1. Introduction.} 
	\medskip
Let $k$ be an algebraically closed field. Let $A$ be a finite-dimensional local $k$-algebra
with radical $J$. The algebra $A$ is said to be {\it short}, provided $J^3 = 0.$
Let $e = \dim J/J^2$ and $s = \dim J^2$. If $A$ is short, then the pair $(e,s)$ is called
the {\it Hilbert type} of $A$. 
The modules to be considered in this paper are usually left $A$-modules, and 
always finitely generated. 
In case $A$ is a short local algebra of Hilbert type $(e,s)$, the local modules 
of length $e$ and Loewy length at most $2$ are of interest;
we call them the {\it atoms} of the short local algebra $A$. 
Several relevant, but standard definitions will be recalled in section 1.5. 

Our aim is to discuss
the existence of non-projective reflexive modules over a short local algebra $A$.
Already in 1971, Ramras [Ra]
asked for a characterization of the rings with the property that
all reflexive modules are projective.

Of course,
if $A$ is self-injective, then all modules are reflexive. Thus, we will assume that $A$
is not self-injective. In cooperation with Zhang Pu, it has been shown that 
the existence of a non-projective reflexive module 
implies that $2 \le s \le e-1$, see [RZ3]: thus  
the dimension of $A$ has to be at least 6, and
if $A$ is 6-dimensional, the Hilbert type of $A$ has to be $(3,2)$. 
Therefore, as a first test case, it seems reasonable to consider in detail 
the short local algebras of Hilbert-type $(3,2)$.
This is done in the present note. From now on, $A$ will be a short local algebra
of Hilbert type $(3,2)$. Let us stress that in our case 
the atoms are the local modules of length 3 and Loewy length 2. 
	\medskip
{\bf 1.1. Theorem.} {\it Let $k$ be an algebraically closed field. Let $A$ be a 
short local $k$-algebra of Hilbert type $(3,2)$ with radical $J$.
The following conditions are equivalent:
\item{\rm(i)} $J^2 = \soc {}_AA = \soc A_A$ and there is no uniform left ideal with length $3.$
\item{\rm(ii)} There is a reflexive atom.
\item{\rm(iii)} There is a non-projective reflexive module.
\par}
	\medskip
The implication (iii) $\implies$ (i) has essentially been shown in [RZ3] (but phrased differently), 
see 7.3 below. Since the implication (ii) $\implies$ (iii) is trivial, 
it remains to show the implication (i) $\implies$ (ii); this is the
target of the present paper. 
The proof of this implication will be based on a study of the 
submodules of ${}_AJ$ of length 2, 3 and 4. 
The proof will be given in 9.6 by showing that there are elements $a\in J\setminus J^2$
such that $A/Aa$ is an atom which is both torsionless and extensionless;
therefore $Aa = \Omega (A/Aa)$ is a reflexive atom
(for any module $M$, we denote by
$\Omega M$ the first syzygy module of $M$, the kernel of a projective cover of $M$).
For the purpose of this paper, an algebra $A$ will be said to be {\it special} 
provided $A$ is a short local algebra of Hilbert type $(3,2)$ satisfying the condition (i).
Thus our aim is to show that {\it any special algebra has reflexive
atoms.} It seems that similar considerations may show that 
any special algebra has even Gorenstein-projective atoms, at least if
$k$ is large, say uncountable. 
	\medskip
{\bf 1.2. The use of $3$-Kronecker modules.} 
The modules of Loewy length at most 2 are just the modules annihilated by $J^2$, thus
the $L(3)$-modules, where $L(3) = A/J^2.$ The $L(3)$-modules correspond nicely to the
$K(3)$-modules, where $K(3)$ is the path algebra of the 3-Kronecker quiver. 
We will recall the definition of the 3-Kronecker quiver 
at the beginning of section 2. The relationship
between $K(3)$-modules and $L(3)$-modules is furnished by
the push-down functor $\pi\:\mod K(3) \to \mod L(3),$ see
the beginning of section 6. 
We say that a $K(3)$-module $W$ with $\bdim W = (3,2)$ 
is {\it special} provided that $W$ is faithful 
and $\End W = k,$ or, equivalently, provided that $W$ is faithful and
has no submodule which is simple injective or which is uniform with length $3$, see 
Proposition 2.3.
Note that {\it a short local algebra $A$ is special iff ${}_AJ = \pi W$
for a special $K(3)$-module $W,$} see 6.5 and 7.3.
	\medskip
{\bf 1.3. Bristles and bars.}
A {\it bristle} (in any length category) is an indecomposable object of length 2.
If the submodule $B$ of the module $M$ is a bristle, we call $B$ a {\it subbristle} of $M$.
If $W$ is a special $K(3)$-module, the number of subbristles of $W$
will be called the {\it bristle type}
of $W$, it turns out that the bristle type of any special $K(3)$-module 
is equal to $1, 2, 3$ or $\infty$,
see 4.2. A 4-dimensional submodule $V$ of a special $K(3)$-module $W$ 
will be called a {\it bar} provided $V$ is not faithful. 
Given a special $K(3)$-module, the subbristles and the bars in $W$ are
the obstacles for atoms to be extensionless or torsionless, see 
Propositions 9.3 and 9.4.
Therefore, we introduce the {\it bristle-bar layout} of $W$ in section 4.1.
As we will see, the shape of the bristle-bar layout only depends on the bristle type of $W$. 
For each bristle type, we will present a coefficient quiver for $W$, 
as well as the corresponding bristle-bar layout in 4.2.
	\medskip
{\bf 1.4. Outline.} Being interested in a special algebra $A$ with radical $J$,
the main object to look at will 
be the projective space $\Bbb PE$, where
$E = J/J^2$. Note that $\Bbb PE$ is a natural parameter space
for the set of isomorphism classes of atoms, see 9.1. 
We will use the fact that $E$ is the top of the left module 
${}_AJ$, as well as 
the top of the right module $J_A$. As we have mentioned, the subbristles and the bars in ${}_AJ$ 
and in $J_A$ have to be studied, thus the left and the right bristle-bar layout of
the algebra $A$.
	
Sections 2, 3, 4, and 5 will be devoted to the 3-Kronecker quiver $K(3)$. 
The relationship between
$K(3)$-modules and $L(3)$-modules is discussed at the beginning of section 6.
The proof of Theorem 1.1 is given at the end of section 9. We should stress that
sections 4 and 5 which provide coefficient quivers for the special $K(3)$-modules
are not essential for the proof of Theorem 1.1.

Section 12 is devoted to commutative algebras where already all 
four bristle types do occur.
And it turns out that for any special commutative algebra, all reflexive atoms are 
Gorenstein-projective.
	\medskip
{\bf 1.5. Some relevant definitions.} For the benefit of the reader, the referees have
suggested to recall relevant definitions (further details may be found in [RZ1] as 
well as in section 2 of 
[RZ3]). Let $M$ be a module.  
The module $M$ is said to be {\it torsionless,} provided it can be embedded into a
projective module. The module $M$ is said to be {\it extensionless,} provided
$\Ext^1(M,{}_AA) = 0.$ A left {\it ${}_AA$-approximation} of $M$ is a map $f\:M \to P$, where
$P$ is projective, such that any map $g\:M \to P'$ with $P'$ projective, can be
factorized as $g = g'f$ with $g':P \to P'$. Such a left ${}_AA$-approximation $f$
is said to be {\it minimal} provided the only direct summand of $P$ containing $f(M)$
is $P$ itself. Any module $M$ has a minimal left 
${}_AA$-approximation $f_M$, and the cokernel
of $f_M$ is denoted by $\mho M$ (and called the {\it agemo} of $M$). 
Of course, $M$ is torsionless iff $f_M$ is 
a monomorphism, and in this case the canonical exact sequence $0 \to M @>f_M>> P \to 
\mho M \to 0$ is said to be an {\it $\mho$-sequence.} The modules of the form $\mho M$ 
are extensionless, and any non-projective indecomposable extensionless module is
of this form. 

The right $A$-module $M^* = \Hom(M, A)$ is called the {\it $A$-dual} of 
$M$. Let $\phi_M\: M \to M^{**}$ be defined by 
$\phi_M(m)(f) = f(m)$ for $m \in M$, $f \in M^*$. Clearly, $M$ is torsionless
provided that $\phi_M$ is injective, and $M$ is called {\it reflexive}
provided that $\phi_M$ is bijective. It is important to be aware 
that $M$ is reflexive iff
both modules $M$ and $\mho M$ are torsionless. 

The module $M$ is said to be {\it Gorenstein-projective} provided there is an
exact complex 
$$
  \cdots @>>> P_1 @>f_1>> P_0 @>f_0>> P_{-1}  @>>> \cdots
$$
with all modules $P_i$ projective, and which remains exact when we form its $A$-dual,
such that $M$ is the image of $f_0$. The Gorenstein projective modules are 
reflexive, but usually reflexive modules will not be Gorenstein-projective. 

Finally, we recall that a module of finite length
is said to be {\it uniform} provided its socle is simple. 
	\bigskip\bigskip
{\bf 2. $K(3)$-modules.}
	\medskip
The {\it $3$-Kronecker quiver} $K(3)$ 
is given by two vertices which we will label $0$ and $1$, and three 
arrows $0 \to 1$. Alternatively, we may start with a 3-dimensional vector space $E$ over $k$ (considered as the arrow space of the quiver), and select, if necessary, a basis of $E$ as
the arrows of the quiver. 
A {\it $3$-Kronecker module} (or $K(3)$-module)
is a representation $M$ of $K(3)$, that means: it 
is of the form $M = (M_0,M_1;\phi_M\:E\otimes_k M_0 \to M_1)$,
where $M_0, M_1$ are (finite-dimensional) vector spaces and $\phi_M$ is a linear map.
In case this map is a canonical one, we just will write $M = (M_0,M_1)$. 
If $M, M'$ are 3-Kronecker modules, a homomorphism $f\:M \to M'$ is a pair $f = (f_0,f_1)$
of linear maps $f_i\:M_i\to M'_i$ (where $i = 0,1$) such that $\phi_{M'}(1\otimes f_0) = f_1\phi_M.$
The $K(3)$-modules form an abelian category denoted by $\mod K(3).$
There are two simple $K(3)$-modules, namely $S_0 = (k,0)$ and $S_1 = (0,k)$.
The module $S_0$ is injective, the module $S_1$ is projective.

If $M$ is a $K(3)$-module, the pair $\bdim M =
(\dim M_0,\dim M_1)$ is called the {\it dimension vector} of $M$
(it is an element of
the Grothendieck group $\Bbb Z^2$ of $\mod K(3)$).
We use the bilinear form $[ -,-]$ on  $\Bbb Z^2$ 
defined by $[ (z_0,z_1),(z_0',z_1')] = z_0z_o'+z_1z_1'-3z_0z_1',$
for $z_0,z_1,z_0',z_1'\in \Bbb Z.$ Note that we have
$[ \bdim M,\bdim M' ] = \dim\Hom(M,M') - \dim \Ext^1(M,M'),$ 
for any pair  $M,M'$ of $K(3)$-modules, see [R1].
	\medskip

{\bf 2.1.} We are going to define some important $K(3)$-modules. 
If $x_1,\dots,x_n$ are elements in a vector space $F$, we denote by 
$\langle x_1,\dots,x_n\rangle$ the subspace of $F$ generated by $x_1,\dots,x_n.$
	\medskip
{\bf (a) Atoms.} For $0\neq a\in E$, let $\widetilde C(a) = \widetilde C(\langle a\rangle) = 
(k,E/\langle a\rangle).$ 
This is an atom, any atom is
obtained in this way, and $\widetilde C(a),\widetilde C(a')$ are isomorphic iff $\langle a\rangle = 
\langle a'\rangle.$ 
{\it The function $a \mapsto \widetilde C(a)$ provides a bijection between $\Bbb PE$
and the set of isomorphism classes of atoms.}

	\smallskip
{\bf (b) Bristles.} If $a_1,a_2$ are linearly independent in $E$, then 
$\widetilde B(a_1,a_2) = \widetilde B(\langle a_1,a_2\rangle ) 
= (k,E/\langle a_1,a_2\rangle).$ This is a
bristle, any bristle is obtained in this way, 
and $\widetilde B(a_1,a_2),\widetilde B(a_1',a_2')$ are isomorphic iff $\langle a_1,a_2\rangle = 
\langle a_1',a_2'\rangle.$ 
{\it The function $\langle a_1,a_2\rangle  \mapsto \widetilde B(a_1,a_2)$
provides a bijection between the set of lines in $\Bbb PE$ and the set of isomorphism classes of 
bristles.}
	\smallskip
{\bf (c) Atoms and bristles.} {\it If $0\neq a\in E$ and $l$ is a line in $\Bbb PE$, then $\langle a \rangle \in l$ iff $\widetilde C(a)$
 maps onto $\widetilde B(l).$}
	\medskip
{\bf 2.2. $K(3)$-modules which are not faithful.}
If $z$ is a non-zero element of $E$, let $\mod K(3)/\langle z \rangle$ be the 
category of all $K(3)$-modules
annihilated by $z$. Then, of course, 
$\mod K(3)/\langle z \rangle$ is equivalent to the category of 2-Kronecker modules.
We will use the well-known structure of the category of 2-Kronecker modules. In particular,
we need the following facts:
	\medskip
{\it The $2$-Kronecker modules with dimension vector $(2,2)$ without a simple direct summand
are either decomposable, then a direct sum of 2 bristles, or else
indecomposable, then an extension of a bristle by itself. An indecomposable 
2-Kronecker-module has a unique
factor module which is a bristle, and this bristle is also the unique subbristle.}
	\medskip
{\it There is only one isomorphism class of indecomposable $2$-Kronecker modules 
with dimension vector $(3,2)$. If $W$ is such a module, any bristle $2$-Kronecker module
occurs as a submodule of $W$.} 
	\medskip
{\it Let $W$ be a $K(3)$-module. There is an atom which 
generates $W$ iff $W$ is not faithful.} 
And if $W$ is a faithful $K(3)/\langle a \rangle$-module, 
then $\widetilde C(a)$ is the only
atom which generates $W$. 
	\medskip
{\bf 2.3. Proposition.} {\it Let $W$ be a $K(3)$-module $W$ with $\bdim W = (3,2)$.
The following conditions are equivalent:
\item{\rm (i)} $\End W = k$.
\item{\rm (ii)} $W$ is indecomposable and has no uniform submodule of length $3$.
\item{\rm (iii)} $W$ has no submodule which is simple injective 
or which is uniform and has length $3$.\par}
	\medskip
Proof. (i) $\implies$ (ii). We
assume that $\End W = k.$ Then $W$ is indecomposable. 
Assume that $U$ is a uniform
submodule of $W$ of length 3. We have $\bdim U = (2,1)$ and therefore
$\bdim W/U = (1,1).$ We use the homological bilinear form $[ -,-]$.
We have $[\bdim W/U,\bdim U] =
[(1,1),(2,1)] =  0$. Since $W$ is indecomposable, $W$ is a non-trivial
extension of $U$ by $W/U$, therefore $\Ext^1(W/U,U) \neq 0.$ It follows that 
$\Hom(W/U,U) \neq 0$, therefore $\End W\neq k$. 

(ii) $\implies$ (iii) is trivial, since a simple injective submodule would be a direct summand.

(iii) $\implies$ (i). We assume that $W$ has no submodule which is simple and injective or which is uniform and has length 3. 
Since $W$ has no direct summand $S_0$, any proper direct decomposition has to be
of the form $W = W'\oplus W''$ with $\bdim W' = (1,1)$ and $\bdim W'' = (2,1).$ However,
in this case, $W''$ is uniform of length 3. Thus, we see that $W$ is indecomposable.
Assume that $\End W \neq k.$ Since $W$ is indecomposable, $\End W$ is a local ring. 
Since it operates on the vector spaces $W_0$ (of dimension 3) and $W_1$ (of dimension 2),
we see that $\End W$ cannot be a division ring. It follows that there is a 
non-zero endomorphism
$\phi$ of $W$ with $\phi^2 = 0.$ Let $I$ be the image of $\phi$ and $K$ the kernel of $\phi$,
thus we deal with submodules $I \subseteq K$ of $W$. Since $I$ is isomorphic to $W/K$, we must
have $\bdim I = (1,1)$ and therefore $\bdim K = (2,1)$. If $K$ is decomposable, then 
$K$ has a direct summand of the form $S_0$, but $S_0$ is injective and thus also $W$ has 
a simple injective direct summand. If $K$ is indecomposable, then $K$ is a 
uniform module of length 3. Both cases are impossible. 
$\s$
	\medskip
We will be interested in the faithful 3-Kronecker modules $W$ 
with $\bdim W = (3,2)$ and $\End W = k$; as we have mentioned, we call 
them just the {\it special} $K(3)$-modules.
	\medskip
{\bf 2.4. Lemma.} {\it Let $W$ be a $K(3)$-module with $\bdim W = (3,2)$ and 
$\End W = k$. Then the indecomposable submodules of $W$ of length $3$ are atoms.}
	\medskip
Proof. This follows directly from 2.3. Namely, 
an indecomposable module of length  3 is either uniform or an atom, but $W$ has no uniform submodule of length 3.
$\s$
	\medskip
If $M, N$ are modules, the largest submodule of $M$ generated by $N$ (thus the sum of the
images of all homomorphisms from $N$ to $M$) is called the {\it trace of $N$ in $M$.}
	\medskip 
{\bf 2.5. Proposition.}
{\it Let $W$ be a special $K(3)$-module. Let $0 \neq a \in E.$ The following
conditions are equivalent:

\item{\rm(i)} The atom $\widetilde C(a)$ has no factor module which is a subbristle of $W$.
\item{\rm(ii)} The trace of $\widetilde C(a)$ in $W$ is an atom.
\item{\rm(ii$'$)} The trace of $\widetilde C(a)$ in $W$ has length $3$.\par}
	\medskip
Let us add: {\it If the trace of $\widetilde C(a)$ in $W$ is an atom, then $\widetilde C(a)$ is
cogenerated by $W$.} The converse is (of course) not true.
	\medskip
Proof. 
(i) implies (ii). We assume that $\widetilde C(a)$ has no factor module which is a subbristle of $W$.
Using the bilinear form $[-,-]$, we see that there is a non-zero homomorphism
$\phi\:\widetilde C(a) \to W$. The image cannot be of length 2, since
otherwise the image would be a subbristle of $W$. Thus $\phi$ is a monomorphism. 
Let $D$ be the image of $\phi$. We claim that $D$ is the trace of $\widetilde C(a)$ in $W.$ 
Otherwise, there is a homomorphism $\psi\:\widetilde C(a) \to W$ with $\psi(a) \notin D.$
Then $N = \Im \phi+\Im\psi$ has length 4, and is generated by $\widetilde C(a)$, thus not faithful.
It follows that $N$ generates a bristle $B$ which is also a subbristle of $N$, see 2.2.
In particular, $B$ is a subbristle of $W$. 
Now $\widetilde C(a)$ generates $N$ and $N$ generates $B$, thus $\widetilde C(a)$ generates $B$.
But this means that $B$ is a factor module of
$\widetilde C(a)$, in contrast to our assumption. 
	\smallskip 
(ii$'$) implies (i). We assume that the trace $D$ of $\widetilde C(a)$ 
in $W$ has length 3. We claim that $D$ is indecomposable. Otherwise, $D$ would have a direct
summand of the form $S_0$ or $S_1$. But $W$ has no submodule of the form $S_0$, since $S_0$
is injective and $W$ is special, and $S_1$ is not generated by $\widetilde C(a).$ Since
$D$ is indecomposable and has length 3, Lemma 2.4 asserts that $D$ is an atom. 
Assume that $\widetilde C(a)$ has a factor module $B$ which is a subbristle of $W.$
Since $D$ is an atom, $B \not\subseteq D$. 
But this contradicts the assumption that $D$ is the trace of $\widetilde C(a)$ in
$W.$
$\s$
	\bigskip\bigskip
{\bf 3. Bristle and bar submodules of special $K(3)$-modules.}
	\medskip

Let $W$ be a  $K(3)$-module with $\bdim W = (3,2).$ We recall that
an indecomposable submodule of $W$ of length $2$ is called a subbristle of $W$. 
A {\it bar} in $W$ is a submodule of $W$ with dimension vector $(2,2)$
which is not faithful. 
	\medskip
A detailed study of the subbristles and bars of a special $K(3)$-module 
will be given in sections 4 and 5. Section 3 provides only some partial information
which will be needed in the proof of Theorem 1.1, namely Propositions 3.7 and 3.8.
	\medskip
{\bf 3.1. The structure of a bar.} {\it Let $W$ be a special $K(3)$-module. 
	\smallskip
\item{\rm (a)} 
 There are three different kinds of bars $V$ in $W$:
\item{$\bullet$} $V$ is the direct sum of two isomorphic bristles, say of bristles isomorphic to $B$ 
 (then all $3$-dimensional submodules are isomorphic to $B\oplus S_0$). 
\item{$\bullet$} $V$ is the direct sum of two non-isomorphic bristles, say $V = B_1\oplus B_2$
 (then $V$ contains no further bristle).
\item{$\bullet$} $V$ is indecomposable and contains precisely one bristle.\par
	\smallskip
\item{\rm (b)} Any bar contains a bristle. 
	\smallskip
\item{\rm (c)} Any bar has a factor module which is a subbristle of $V,$ thus of $W$.}
	\smallskip
\item{\rm (d)} {\it 
All indecomposable submodules of length $3$ of a bar are isomorphic.}
(Note that indecomposable submodules of length 3 do not exist iff $V$ is the direct sum of two isomorphic bristles.) \par
	\medskip
Proof. This follows directly from the classification of the $2$-Kronecker modules, as mentioned
in 2.2.
Namely, assume that the bar $V$ is annihilated by $0\neq z \in E$, thus it may be considered as
a $K(3)/\langle z \rangle$-module (thus as a $2$-Kronecker module).
In this way, $V$ is just a 4-dimensional $2$-Kronecker module without a simple
direct summand. For the proof of (d), let $U$ be an indecomposable submodule of $V$ of length 3.
Now any submodule of $V$ is again a
$K(3)/\langle z \rangle$-module. Since $W$ has no uniform submodule of length 3,
we see that $U$ has to be local.  Thus $U$ is a
projective cover of $S_0$ considered as a $K(3)/\langle z \rangle$-module.
$\s$
	\medskip
{\bf 3.2. Lemma.}
{\it Let $W$ be a special $K(3)$-modules and $B, B'$ different subbristles of $W$,
then $B\cap B' = 0$, and $B\oplus B'$ is a bar.}
	\medskip
Proof: If $B\cap B' \neq 0,$ then $B+B'$ is a submodule of $W$ with dimension vector $(2,1)$,
thus either uniform or isomorphic to $B\oplus S_1$, impossible. 
Thus $B\cap B' = 0$ and therefore
$B + B' = B\oplus B'$. The annihilator of a bristle is a 2-dimensional subspace of $E$.
The annihilator of $B\oplus B'$ is the intersection of the annihilator of $B$ and the
annihilator of $B'$, thus non-zero.
$\s$
	\medskip 
{\bf 3.3. Lemma.} {\it Let $W$ be a special $K(3)$-module. Let $N, N'$ be different bars
in $W$. Then there is a bristle $B$ such that $N\cap N' = B+\soc W.$}
	\medskip
Proof. Both $N, N'$ contain $\soc W$. Thus, the dimension vector of $N\cap N'$ is $(1,2).$
Assume that $N\cap N'$ is indecomposable, say isomorphic to $(k,E/\langle z \rangle).$ Since
$N$ is not faithful, and has a submodule of the form $(k,E/\langle z \rangle),$ it follows
that $N$ is annihilated by $z$. Similarly, $N'$ is annihilated by $z$. Since $W = N+N'$, we
see that $W$ is annihilated by $z$. This contradicts the assumption that $W$ is faithful.

Since $N\cap N'$ is decomposable, and $W$ has no submodule isomorphic to $S_0$, it follows
that $N\cap N'$ contains an indecomposable module of length $2$, a bristle.
$\s$
	\medskip 

{\bf 3.4. Lemma.} {\it Let $W$ be a special $K(3)$-module.
If there are two isomorphic subbristles $B, B'$, then $B\cap B' = 0$ and 
$B\oplus B'$ is the only bar in $W$.
The indecomposable direct summands of $B\oplus B'$
are all the subbristles of $W$, and all subbristles of $W$ are isomorphic to $B$.} 
	\medskip
Proof: According to 3.2, we have $B \cap B' = 0$ and $V = B_1\oplus B_2$ is a bar.
Since $B$ and $B'$ are isomorphic, all submodules of 
$V$ of length 3 are of the form $B\oplus S_1.$ 

Assume that there is a bristle $B''$ outside of $V$. Then $W = V+B''$ is annihilated by
$\Ann(B)\cap \Ann(B''),$ and 
$\Ann(B)\cap \Ann(B'') \neq 0$, thus $W$ is not faithful - a contradiction. 
This shows that the indecomposable direct summands of $V$ are the
subbristles of $W$. In particular, all subbristles of $W$ are isomorphic to $B$.

Now, assume that there is a second bar $V'$. According to 3.3, 
$V\cap V' = B''+\soc W$ for some bristle $B''$. Since $B''$ is a subbristle of $W$,
we see that $B''$ is isomorphic to $B$. Now $V'$ is not faithful, say annihilated by $0\neq z \in E.$
Then $z$ annihilates $B''$, thus $B$ and $B'$. Since $W = V+V' = B+B'+V'$, we see that
$z$ annihilates $W$, impossible.
$\s$
	\medskip
{\bf 3.5. Lemma.} {\it Let $W$ be a special $K(3)$-module.
If there are (at least) three bristles $B_1, B_2, B_3$, not all isomorphic,
then they are pairwise non-isomorphic, they are the only bristles, 
the modules $B_i\oplus B_j$ with $i\neq j$ are the only bars and $W =B_1+B_2+B_3$.}
	\medskip
Proof: The bristles are pairwise non-isomorphic, according to 3.4. If $B_1\cap B_2 \neq 0,$
then $B_1+B_2$ is indecomposable, thus uniform. But this is impossible. 
Therefore $B_1\cap B_2 = 0,$ thus $B_1+B_2 = B_1\oplus B_2.$ Similarly, we have 
$B_1+B_3 = B_1\oplus B_3$ and $B_1+B_3 = B_2\oplus B_3.$ 
Since $B_1, B_2$ are non-isomorphic, $B_1$ and $B_2$ are the only subbristles of
$B_1\oplus B_2$, therefore $B_3 \not\subseteq (B_1+B_2).$ It follows that 
$W = B_1+B_2+B_3.$ 

Assume that $V$ is a bar. Then $V\cap (B_1\oplus B_2)$ is of the form $B+\soc W$ with
$B$ a bristle, say, without loss of generality $B = B_1$.
Also $V\cap B_2\oplus B_3$ is of the form $B'+\soc W$ for some bristle $B'$, thus 
$B'$ is either $B_2$ or $B_3$. In the first
case, we have $V = B_1\oplus B_2,$ in the second case, we have $V = B_1\oplus B_3.$
Thus, there are just three bars $B_1+B_2,\ B_1+B_3,\ B_2+B_3.$

Finally, consider any subbristles $B$. If $B \neq B_1,$ then 
$B + B_1$ is a bar, thus one of the three bars $B_1+B_2,\ B_1+B_3,\ B_2+B_3.$
The only subbristles of these bars are the bristles $B_1,B_2,B_3$, thus $B$ 
is one of them.
$\s$
	\medskip
{\bf 3.6. Lemma.} {\it A special $K(3)$-module $W$ has at most one 
indecomposable bar submodule.}
	\medskip
Proof. Let $W$ be a special $K(3)$-module. 
Assume $N\neq N'$ are indecomposable bars in $W$. 
According to 3.3, $N\cap N' = B+\soc W$ for some bristle $B$. 

Since $B$ is a bristle, and $N$ is an indecomposable 2-Kronecker module, the classification of the indecomposable 2-Kronecker modules shows that $N/B \simeq B.$ Similarly, $N'/B \simeq B.$
We claim that $W/B$ has a direct summand isomorphic to $S_0.$ Namely, we have
$W/B = (N+N')/B = N/B + N'/B.$ Since $N/B$ and $N/B'$ are isomorphic, $N/B$ is a direct summand
of $N/B+N'/B$ and a direct complement has to be of the form $S_0.$ Thus $W/B = N/B \oplus 
U/B$ for some submodule $U$ with $B \subseteq U$ and $U/B \simeq S_0.$ If $U$ is decomposable,
then $U \simeq B\oplus S_0$, thus $S_0$ is a submodule of $W$. Or else $U$ is uniform.
Both cases are impossible, according to 2.4.
$\s$
	\medskip
{\bf 3.7. Proposition.} {\it A special $K(3)$-module $W$ has at most three 
isomorphism classes of subbristles and at most three bars.}
	\medskip
Proof. If $W$ has subbristles $B \neq B'$ which are isomorphic, then we apply 3.4
in order to see that $W$ has only one isomorphism class of subbristles and only one bar. 
If $W$ has at least 3 subbristles which are pairwise non-isomorphic, then we apply 3.5
in order to see that $W$ has only 3 subbristles (thus 3 isomorphism classes of subbristles)
and only 3 bars. 

Thus, we can assume that $W$ has at most 2 subbristles. 
We claim that $W$ has at most two bars.
If $W$ has two subbristles $B, B'$, then $B\oplus B'$ is a bar, and this is 
the only decomposable bar. If $W$ has just one subbristle, then there cannot be any decomposable bar. 
This shows that if $W$ has at most 2 subbristles, there is at most one decomposable bar. 
And according to 3.6, there is at most one indecomposable bar. 
$\s$
	\medskip
{\bf Remark.} As we will see later, the number of 
isomorphism classes of subbristles of a special 3-Kronecker module is the same
as the number of bar submodules, see 4.6.
	\medskip
{\bf 3.8. Proposition.} {\it Let $M$ be a $K(3)$-module with $\bdim M = (3,2).$
Any subbristle of $M$ is contained in a bar submodule.}
	\medskip
Proof. Let $B$ be a subbristle of $M$, annihilated by linearly independent elements $y,z \in E$.
We choose elements $c,d$ of $M$ such that $B, c, d$ generate $M$. Let 
$F = \langle y,z\rangle$ and $T = \langle c,d\rangle.$ 
We consider the multiplication map $\mu\:F\otimes T \to \soc M.$
This is a linear map, say with kernel $K$ of dimension at least $2$
(since $F\otimes T$ has dimension 4). 
The set of decomposable tensors in $F\otimes T$ is a subvariety $\Cal Y$ 
isomorphic to $\Bbb P^1\times\Bbb P^1,$ thus $\dim \Cal Y = 2.$
Let $\Cal K = \Bbb PK$, it has dimension 1. 
Thus $\Cal K$ and $\Cal Y$ are subvarieties of $\Bbb P(F\otimes T)$ which intersect non-trivially, 
since $\dim \Cal K + \dim \Cal Y \ge \dim \Bbb P(F\otimes T)$.
Thus there are elements $0 \neq e \in F$ and $0\neq t \in T$ with $\mu(e\otimes t) = 0.$ Let
$V = B + At + \soc W$. Then $V$ is 4-dimensional, contains $B$ and is 
annihilated by $e$ (note that
$e\in F$ annihilates not only $t$ but also $B$, and of course $\soc W$). This shows that $V$ is a bar.
$\s$
	\medskip
{\bf 3.9. Addendum.} {\it  Any special $K(3)$-module contains a bristle.} 
	\medskip
Remark. An immediate consequence of 3.9 and 3.8 is: 
{\it Any special $K(3)$-module contains a bar.} This is the essential new information,
because the known classification of the $2$-Kronecker modules shows 
that any bar contains a bristle.
	\medskip
Proof of 3.9. It has been shown in [R3] that any
$K(3)$-module $M$ with $\bdim M = (2,3)$ has a submodule $U$ with dimension vector $(1,2)$,
thus $M$ has a factor module with dimension vector $(1,1).$ Using duality, we see that 
any $K(3)$-module with
dimension vector $(3,2)$ has a submodule with dimension vector $(1,1).$
      
For the convenience of the reader, let us copy the proof of [R3]. 
Let $M = (M_1,M_2)$ be a $K(3)$-module with $\bdim (2,3)$. We first 
show that there are non-zero elements 
$m\in M_1$ and $\alpha\in E $ such that $\alpha m = 0$. 
Let $K$ be the kernel of the map $f\:E  \otimes_k M_1 \to M_2$.
We have $\dim E \otimes_k M_1 = 3\cdot 2 = 6$. Since $\dim M_2 = 3,$ 
it follows that $\dim K \ge 3.$ The projective
space $\Bbb P(E \otimes M_1)$ has dimension $5$, the decomposable tensors in $E \otimes M_1$
form a closed subvariety $\Cal V$ of $\Bbb P(E \otimes M_1)$ of dimension $(3-1)+(2-1) = 3$.
On the other hand, $\Cal K = \Bbb P(K)$ is a closed subspace of $\Bbb P(E \otimes M_1)$ 
of dimension $2$. It follows that 
$$
 \dim(\Cal V \cap \Cal K) \ge 3+2-5 \ge 0,
$$
thus
$\Cal V \cap \Cal K$ is non-empty. As a consequence, we get non-zero 
elements $m\in M_1, \alpha\in E $ such that $\alpha m = 0,$ as required.

Given non-zero elements $m\in M_1$ and $\alpha\in E $ such that $\alpha m = 0$, 
the element $m$ generates a submodule $U'$ which is annihilated by $\alpha$, thus
$\bdim U' = (1,u)$ with $0\le u \le 2$.
Since $\dim M_2 \ge 3$, there is a semi-simple submodule $U''$ 
of $M$ with dimension vector $(0,2-u)$ such that
$U'\cap U'' = 0$. Let $U = U'\oplus U''$. This is a submodule of $M$ with dimension vector
$\bdim U = \bdim U'\oplus U'' = (1,2)$.
$\s$
	\medskip
{\bf 3.10.} 
As we have mentioned already in the introduction, the number of subbristles of a special
$K(3)$-module $W$ will be called the bristle type of $W$.
	\bigskip\bigskip
{\bf 4. Coefficient quivers and bristle-bar layouts.}
	\medskip
We are going to describe coefficient quivers for the special $K(3)$-modules
as well as the corresponding bristle-bar layouts.
	\medskip
{\bf 4.1.} 
The {\it bristle-bar layout} 
of a special 3-Kronecker module $W$ is given by the projective plane
$\Bbb P(W_0)$ with marks which highlight the elements $\langle a \rangle\in \Bbb P W_0$ with $Aa$ being a subbristle of $W$
(say using squares $\ssize \blacksquare$)
as well the top lines $t(N) = \{\langle a \rangle\in \Bbb P W_0\mid Aa \subset N\}$, 
where $N$ is a bar in $W$ (say by drawing 
$t(N)$  as a solid line). 
To any bar $N$, we add the annihilator $\langle z_1,\dots, z_t\rangle$ of $N$ in $E$ (here,
$z_1,\dots,z_t$ are elements of $E$). (Note that this kind of visualization was used
already in our joint paper [RZ2] with Zhang Pu.) 
	\medskip
{\bf Convention:} If $W_0 = k^3,$ we draw the plane $\Bbb PW_0$ as follows: the
left corner is $(100)$, the right corner is $(010)$, thus the upper corner is $(001).$
$$
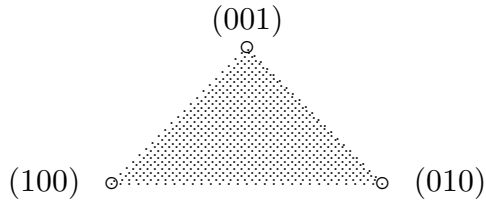
{\beginpicture
    \setcoordinatesystem units <1.5cm,1.5cm>
\setdots <1mm>
\plot 0 0  2 0  1 1  0 0 /
\put{$(100)$} at -.5 0
\put{$(010)$} at 2.5 0
\put{$(001)$} at 1 1.2 
\multiput{$\circ$} at 0 0  2 0  1 1 /
  \setshadegrid span <.4mm>
  \vshade 0 0 0 <,z,,> 1 0 1 <z,,,> 2 0 0 /
\endpicture}
$$
	\medskip
{\bf 4.2. The special $K(3)$-modules.} Recall that the bristle type of 
of a special $K(3)$-module $W$ is the number of subbristles of $W$.
Given a special $K(3)$-module $W$, we are going to write down a coefficient quiver, 
as defined in [R2]. 
	\medskip
{\bf Proposition.} {\it The bristle type of a special 
$K(3)$-module $W$ is $1, 2, 3$ or $\infty.$
Depending on its bristle type, $W$ can be described by one of the following coefficient quivers.}
On the right, we show
the corresponding bristle-bar layout.
$$
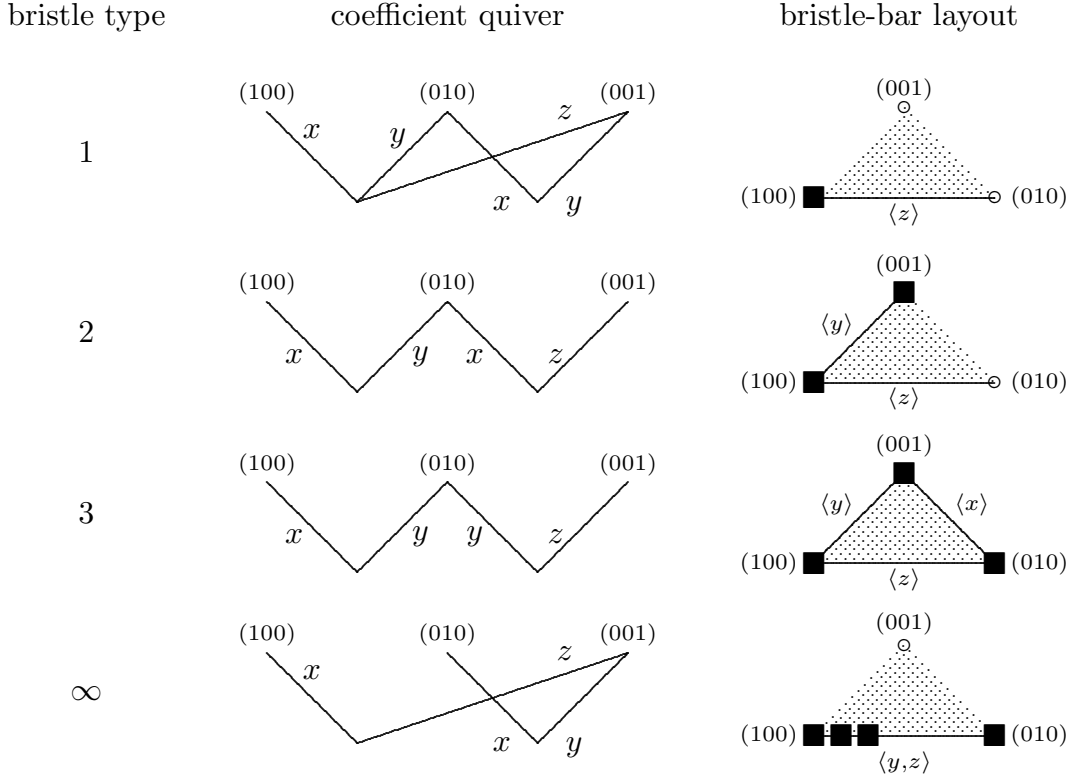
{\beginpicture
    \setcoordinatesystem units <1cm,1cm>
\put{bristle type} at -1 5.5
\put{coefficient quiver} at 3 5.5
\put{bristle-bar layout} at 8 5.5
\put{
\beginpicture
  \multiput{$\blacksquare$} at 0 0  /
  \multiput{} at 2 0 /
  \plot 0 0  2 0  /
  \multiput{$\circ$} at 2 0   1 1 /
  \setdots <1mm> 
  \plot 0 0   2 0  1 1 /
  \put{$\ssize (100)$} at -0.5 0
  \put{$\ssize (001)$} at 1 1.2
  \put{$\ssize (010)$} at 2.5 0
  \setshadegrid span <.5mm>
  \vshade 0 0 0 <,z,,> 1 0 1 <z,,,> 2 0 0 /
\put{$\ssize \langle z \rangle$} at 1 -.18
\endpicture} at 8 4
\put{\beginpicture
  \multiput{} at  0 0  4 1 /
  \plot 0 1  1 0  2 1  3 0  4 1  1 0 /
  \multiput{$x$\strut} at 0.5 0.8  2.6 0 /
  \multiput{$z$} at  3.3 1 /
  \multiput{$y$\strut} at 1.45 0.75  3.4 0 /
  \put{$\ssize (100)$} at 0 1.2
  \put{$\ssize (010)$} at 2 1.2
  \put{$\ssize (001)$} at 4 1.2
\endpicture} at 3 4
\put{1} at -1 4

\put{
\beginpicture
  \multiput{$\blacksquare$} at 0 0   1 1 /
  \put{$\circ$} at 2 0
  \plot 1 1  0 0  2 0    /
  \setdots <1mm> 
  \plot 0 0  1 1  2 0 /
  \put{$\ssize (100)$} at -0.5 0
  \put{$\ssize (001)$} at 1 1.3
  \put{$\ssize (010)$} at 2.5 0
 \setshadegrid span <.5mm>
 \vshade 0 0 0 <,z,,> 1 0 1 <z,,,> 2 0 0 /
\put{$\ssize \langle y \rangle$} at 0.25 0.65
\put{$\ssize \langle z \rangle$} at 1 -.18
\endpicture} at 8 2
\put{\beginpicture
  \multiput{} at  0 0  4 1 /
  \plot 0 1  1 0  2 1  3 0  4 1 /
  \multiput{$x$} at 0.3 0.4  2.3 0.4 /
  \multiput{$y$} at 1.7 0.4 /
  \multiput{$z$} at 3.2 0.4 /
  \put{$\ssize (100)$} at 0 1.2
  \put{$\ssize (010)$} at 2 1.2
  \put{$\ssize (001)$} at 4 1.2
\endpicture} at 3 2
\put{2} at -1 2

\put{
\beginpicture
  \multiput{$\blacksquare$} at 0 0  1 1  2 0 /
  \plot 0 0  2 0   1 1  0 0 /
  \setdots <1mm> 
  \plot 0 0  1 1  2 0 /
  \put{$\ssize (100)$} at -0.5 0
  \put{$\ssize (001)$} at 1 1.3
  \put{$\ssize (010)$} at 2.5 0
 \setshadegrid span <.5mm>
 \vshade 0 0 0 <,z,,> 1 0 1 <z,,,> 2 0 0 /
\put{$\ssize \langle x \rangle$} at 1.75 0.65
\put{$\ssize \langle y \rangle$} at 0.25 0.65
\put{$\ssize \langle z \rangle$} at 1 -.18
\endpicture} at 8 0
\put{\beginpicture
  \multiput{} at  0 0  4 1 /
  \plot 0 1  1 0  2 1  3 0  4 1 /
  \multiput{$x$} at 0.3 0.4 /
  \multiput{$y$} at 1.7 0.4  2.3 0.4 /
  \multiput{$z$} at 3.2 0.4 /
  \put{$\ssize (100)$} at 0 1.2
  \put{$\ssize (010)$} at 2 1.2
  \put{$\ssize (001)$} at 4 1.2
\endpicture} at 3 0
\put{3} at -1 0

\put{
\beginpicture
  \multiput{$\blacksquare$} at 0 0  0.3 0  0.6 0  2 0 /
  \plot 0 0  2 0 /
  \put{$\circ$} at 1 1 
  \setdots <1mm> 
  \plot 0 0  1 1  2 0 /
  \put{$\ssize (100)$} at -0.5 0
  \put{$\ssize (001)$} at 1 1.25
  \put{$\ssize (010)$} at 2.5 0
 \setshadegrid span <.5mm>
 \vshade 0 0 0 <,z,,> 1 0 1 <z,,,> 2 0 0 /

\put{$\ssize \langle y,z \rangle$} at 1 -.3
\endpicture} at 8 -2
\put{\beginpicture
  \multiput{} at  0 0  4 1 /
  \plot 0 1  1 0  4 1  3 0  2 1  /
  \multiput{$x$\strut} at 0.5 0.8  2.6 0 /
  \multiput{$z$} at  3.3 1 /
  \multiput{$y$\strut} at  3.4 0 /

  \put{$\ssize (100)$} at 0 1.2
  \put{$\ssize (010)$} at 2 1.2
  \put{$\ssize (001)$} at 4 1.2
\endpicture} at 3 -2
\put{$\infty$} at -1 -2
\endpicture}
$$
	\medskip
The proof of Proposition 4.2 will be given in section 5. In the present section 4, we 
verify that a Kronecker module with one of the four coefficient quivers presented in the
table is special and has the described bristle-bar layout. We also present 
a bijection between the subspaces $\langle z \rangle$ of $E$, where $\dim zW = 1$
and the subbristles of $W$, see 4.5.
	\medskip
{\bf 4.3.} {\it Let $W$ be a $K(3)$-module 
as exhibited in} 4.2. {\it Let $0 \neq a \in W_0$.
If $\langle a \rangle$ is marked by a black square, then $a$ generates a bristle, otherwise it generates 
an atom.}
	\medskip
Proof. It is straightforward to see
that the bristles displayed in the layout actually do exist. We only write
down how to verify that the remaining non-zero elements $a\in W_0 = k^3$ generate atoms. 
We identify $W_1 = k^2$, by choosing as first basis element of $W_1$
the element $x(100)$, as second one the element $x(010)$ if the bristle type is $1, 2$ or $\infty$
and $y(010)$ if the bristle type is $3$.
	\smallskip
Type 1. 
First, assume that $\langle a \rangle$ lies on the shown bar, but is different from 
$\langle 100\rangle,$
say $a = (\alpha,1,0)$ for some $\alpha\in k.$ Then $xa = (\alpha,1)$ and $ya = (1,0)$,
thus $a$ generates an atom. Second, assume that $\langle a \rangle$ does not lie on the bar,
say $a = (\alpha,\beta,1)$ for some $\alpha,\beta \in k.$ Then
$za = (1,0)$ and $ya = (\beta,1)$. Again, $a$ generates an atom.
	\smallskip
Type 2. 
The elements $(100)$ and $(001)$ generate bristles. Since these bristles are 
non-isomorphic, $(\alpha,0,\gamma)$ with $\alpha\neq 0, \gamma \neq 0$ generates an atom. 
Finally, let $a = (\alpha,1,\gamma).$
Then $x a = (\alpha,1)$ and $z a = (1,0)$ show that $a$ generates an atom.
Thus, there are precisely 2 bristles.
	\smallskip
Type 3. The elements $(100), (010)$ generate pairwise non-isomorphic bristles,
thus the elements $(1,\beta,0)$ with $\beta \neq 0$ generate atoms. Similarly, the
elements $(1,0,\gamma)$ and $(0,1,\gamma)$ with $\gamma\neq 0$ generate atoms.
Finally, consider $a = (\alpha,1,\gamma)$ with $\alpha\neq 0$ and $\gamma\neq 0.$ Then
$xa = (\alpha,0)$ and $za = (0,\gamma)$ shows that $a$ generates an atom.
	\smallskip
Type $\infty$. The elements on the bar generate bristles. Let $a$ be outside the bar, thus
we can assume $a = (\alpha,\beta,1)$. Then $za = (1,0)$ and $ya = (0,1)$, thus $a$
generates a bristle.
$\s$
	\medskip
{\bf 4.4.} {\it Let $W$ be a $K(3)$-module as exhibited in} 4.2. {\it Then $W$
is a special $K(3)$-module.}
	\medskip
Proof. The module $W$ has no submodule which is simple injective, since the non-zero
elements of $W_0$ generate bristles and atoms. The module $W$ has no uniform submodule $V$
of length 3. Namely, $t(V)$ would be a line consisting of pairwise non-isomorphic 
bristles. Lines of bristles do exist only in type $\infty$, but in this case, all the bristles
are isomorphic. 

According to Proposition 2.3, it follows that $\End W = k,$ 
thus, $W$ is indecomposable. As we have mentioned in 2.2, an indecomposable 
3-Kronecker module with dimension vector $(3,2)$ which is not faithful has 
infinitely many pairwise non-isomorphic subbristles. Thus $W$ has to be faithful.
 $\s$.
	\medskip
	
{\bf 4.5. Proposition.} {\it Let $W$ be a $K(3)$-module with $\bdim W = (3,2).$
\item{\rm (a)} If $W$ has no simple injective submodule and 
$z\in E$ satisfies $\dim zW = 1$, then $zW$ is the socle of a subbristle.
\item{\rm (b)} If $W$ is special, then we obtain in this way a bijection 
between the subspaces $\langle z\rangle$ 
generated by the elements $z\in E$ with $\dim zW = 1$ and the subbristles of $W$.\par}
	\medskip
Proof. (a) Let $z\in E$ with $\dim zW = 1.$ Of course, $U = (0,zW)$ is a submodule of $W$
and $W/U$ is a 3-Kronecker module annihilated by $z$, thus a 2-Kronecker module for
$E/\langle z \rangle$.
The dimension vector of $W/U$ is $(3,1)$. It must have a direct summand isomorphic to $S_0$. Thus there
is a submodule $V \subseteq W$ with $U\subseteq V$ and $V/U$ isomorphic to $S_0.$
Of course, $\bdim V = (1,1).$ If $V$ would be decomposable, $S_0$ would be a submodule of $V$, thus
of $W$. But this is impossible, since we assume that $W$ has no simple injective submodule.
Therefore $U$ is indecomposable,
thus a bristle. The socle of $U$ is just $zW.$
	\smallskip
(b) We assume now that $W$ is special. First, let us start with a subbristle $B$ of $W$. 
Let $U = \soc B$. Then $W/U$ has $B/U$ as a submodule.
Since $W/U$ is isomorphic to $S_0$, we see that $B/U$ is a direct summand, thus there is 
$V\subseteq W$ with $U \subseteq V$ and $W/U = B/U\oplus V/U.$ The module $V$ has $\bdim V = (2,1)$,
thus it cannot be faithful: there is $0\neq z \in E$ with $zV = 0.$ But this means that 
$zW \subseteq U.$ Since $W$ is special, $W$ is faithful, thus
we must have $zW \neq 0$, therefore $zW = U.$

It remains to be seen that given a bristle $B$ and elements $z,z'$ with $zW = z'W = \soc B$,
the subspaces generated by $z$ and by $z'$ are equal. 
Assume not. Then there is an element $x\in E$
such that $x,z,z'$ is a basis of $E$. Let $U = zW = z'W.$ 
Since $W_0$ is 3-dimensional, there is a 2-dimensional subspace $X$ of $W_0$ with 
$xX \subseteq U$. We claim that $(X,U)$ is a submodule of $W$. Namely, since $U = zW = z'W$,
we have $zX \subseteq U$ and $z'X \subseteq U$. By assumption, we also have $xX \subseteq U.$
Note that $\bdim (X,U) = (2,1).$ 
Since $W$ is special, $W$ has no submodule with dimension vector $(2,1)$.
$\s$
	\medskip
{\bf 4.6. Corollary.}
{\it Let $W$ be a special $K(3)$-module. 
The number $\theta$ of bristles in $W$ 
is the same as the number of subspaces $\langle z \rangle$ of $E$
with $\dim zW = 1$, namely $1,2,3$ or $\infty$.}  
Recall that this number $\theta$ is called the bristle type of $W$.  
{\it The number of bars in $W$ is less or equal to $\theta$ (the number of subbristles of $W$).}
	\medskip
Actually, when we have finished the proof of 4.2, we will know that 
{\it the number of bars in a special Kronecker module $W$ is equal to the number of 
{\bf isomorphism classes} of subbristle of $W$.}
	\medskip
Proof of Corollary. The first assertion is just 4.6 (b). For the second assertion,
one only has to observe: For any $z\in E$ with $\dim zW = 1$, the submodule $(\Ker z,W_0)$
is a bar in $W$. 
$\s$
	\medskip
{\bf 4.7. Remark.} 
We have seen in 4.5 that given a special 3-Kronecker module $W$ and $z\in E$ with $\dim zW = 1$,
then $zW = \soc B$ for some bristle $B$. Let us stress that there are two different possibilities,
namely, we may have either $zB = 0$ or else $zB \neq 0$ (and then $zB = \soc B)$.

Both possibilities occur, as the special modules of bristle type 2 
(exhibited in 4.2) show:
There are two bristles, namely $B$ generated by $(100)$ and $B'$ generated by $(001)$.
The element $y\in E$ satisfies $yW = \soc B$ and we have $yB = 0.$
The element $z\in E$ satisfies $zW = \soc B'$ and we have $zB' \neq 0.$
	\medskip
{\bf 4.8.} {\it Let $W$ be a special $K(3)$-module as exhibited in} 4.2. {\it Then
all bars in $W$ are shown in the corresponding layout on the right.}
	\medskip
Proof. Again, it is clear that the lines marked in the layout are bars. We have to show
that there are no additional ones.
	\smallskip
Type $\theta = 1, 2, 3.$ Here we have $\theta$ bristles. According to 4.4, $W$ is a special 
$K(3)$-module. Thus, according to 4.6, there are at most $t$ bars. Since $t$ lines
are marked as bars, we see that these are the only bars. 
	\smallskip
Type $\infty$. Since the submodules generated by $(100)$ and $(010)$ are isomorphic bristles,
Lemma 3.4 asserts that $W$ has only one bar.
$\s$

	\medskip
{\bf 4.9. Remark.} Let $W$ be a 3-Kronecker module as exhibited in 3.2.
We see: if the number of subbristles is finite, then 
the number of subbristles is the same as the number
of bars, namely $1,$ $2,$ or $3$ and that any bar contains a bristle. 
Also, if a bar contains at least $3$ bristles, then it contains infinitely many bristles
and there is no other bar. Finally, a line with at least two bristles is a bar. 

Altogether, these assertions imply:
{\it If there are precisely three bristles, the bars form a triangle} (they intersect 
in pairwise different vertices).
	\medskip
Proof. Namely, the three bristles cannot lie on a single bar (since otherwise there would exist
infinitely many bristles), thus they are the corners of a triangle. This yields 
3 bars, so there cannot be a further bar. $\s$
	\bigskip\bigskip

{\bf 5. Proof of Proposition 4.2.}
	\medskip
{\bf 5.1. Type 1: Precisely one bristle.} 
	\medskip 
Proof. Let $B \subset N$, where $B$ is the bristle and $N$ a bar. Since there is only one
bristle, $N$ is indecomposable. Let $a$ be a generator of $B$
and assume that $z$ generates the annihilator of $N$. 
According to 4.5 (a), we have $zW = \soc B$. Since $zN = 0$, we have $za = 0.$
Since $B$ is a bristle and annihilated by $z$,  there is 
a basis $x, y, z$ such that also $yB = 0.$
Now $a$ is not in the kernel of $x$, thus $d = xa$ is a non-zero element in $W_1.$
Since $N$ is an indecomposable bar with subbristle $B$, also $N/B$ is isomorphic to $B$,
thus $x$ maps $N\cap W_0$ bijectively onto $W_1.$ 

Since $x$ has rank 2, its kernel has dimension 1. Let $c$ be a non-zero element
in the kernel of $x$.  
According to 4.5 (a), the image of $z$ coincides with the socle of $B$.
Thus, replacing, if necessary, $c$ by a non-zero multiple, we can assume that $zc = d.$ 
Note that $d = zc$ and $e = yc$ have to be linearly independent in $W_1$, since otherwise $c$
would generate a second bristle. 

Choose $b'\in N\cap W_0$ such that $yb' = d$ (since $ya = 0$, the elements $a, b'$ are
linearly independent). Write $xb'$ as a linear combination of $d$ and $e$, say
 $xb' = \delta d + \epsilon e.$ Note that $\epsilon \neq 0$, since
otherwise $b'$ would generate a bristle. Let $b = -\delta a+b'.$
Then 
$$ 
   xb =  -\delta xa +   xb' = 
       -\delta d    +   \delta d + \epsilon e = \epsilon e.
$$
and
$$
   yb =  -\delta ya + yb' = d. 
$$
We obtain the following coefficient quiver
$$
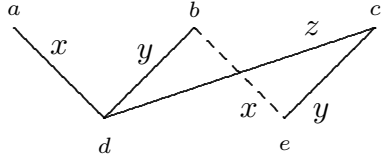
{\beginpicture
    \setcoordinatesystem units <1cm,1cm>
  \multiput{} at  0 0  4 1 /
  \plot 0 1  1 0  2 1 /
  \plot 3 0  4 1  1 0 /
  \setdashes <1mm>
  \plot 2 1  3 0 /
  \multiput{$x$\strut} at 0.5 0.8  2.6 0.1 /
  \multiput{$z$} at  3.3 1 /
  \multiput{$y$\strut} at 1.45 0.75  3.4 0.1 /
  
  \put{$\ssize a$} at 0 1.2
  \put{$\ssize b$} at 2 1.2
  \put{$\ssize c$} at 4 1.2
  \put{$\ssize  d$} at 1 -.3
  \put{$\ssize  e$} at 3 -.3
\endpicture} 
$$
where the dashes line indicates that here a non-zero scalar is involved, 
namely $xb = \epsilon e.$ In order to
remove the scalar, we replace $x$ by $\epsilon^{-1} x$ and $a$ by $\epsilon a.$ Then we get the 
same picture, with the dashed line replaced by a solid one. 
\vglue-.2cm
$\s$
	\medskip
{\bf 5.2. Type 2. Precisely two bristles.} 
	\medskip
Proof. Let $B \neq B'$ be the bristles of $W$. 
According to 3.4, they are not isomorphic and $B\cap B' = 0.$ 
The annihilator of $B\oplus B'$ is $\langle y \rangle$ for some
non-zero element $y$, therefore $\dim yW = 1.$
According to 4.6, there is an element $z\in E$ such that also $\dim zW = 1$, with 
$y,z$ being linearly independent.
The submodule $\Ker z$ is a bar. 
Note that $\Ker y \neq \Ker z$, since otherwise $B\oplus B'$ 
would be annihilated by $y$ and $z$,
and then $B$ and $B'$ would be isomorphic bristles. 
According to 3.3, the intersection $\Ker z\cap \Ker y$ contains a bristle. 
Without loss of generality, $B \subseteq \Ker z\cap \Ker y$,
and then $B'$ is not contained in $\Ker z.$

Since $B'$ is not contained in $\Ker z$, we see that 
$zB'$ is the socle of $B',$ therefore $zW =  \soc B'.$ 
On the other hand, using again 4.6, $yW$ is the socle of a bristle different from $B'$,
thus $yW = \soc B.$

Choose $x\in E$ such that $x,y$
are linearly independent and annihilate $B'$. Since $z$ does not annihilate $B'$, we
see that $x,y,z$ is a basis of $E$. We know that $yW$ and $zW$ are both one-dimensional.
Since there are only 2 bristles, 4.6 asserts that we have $\dim xW = 2,$ thus $xW = W_1.$ 

The bristle $B$ is annihilated by $y,z$, the bristle $B'$ is annihilated by $x,y$.
Since these are the only bristles, no non-zero element in $W_0$ is annihilated by $x$ and $z$.
Thus $(\Ker z)\cap W_0$ is mapped under $x$ bijectively to $W_1.$ 

Let $c$ be a generator of $B'$ and $e = zc.$
Choose $b \in (\Ker z)\cap W_0$ such that $xb = e.$
Let $d = yb.$ Then $d$ is a non-zero element of the socle of $B$. Let $a$ be a generator of $B$
with $xa = d.$ Altogether we have obtained the following coefficient quiver:
$$
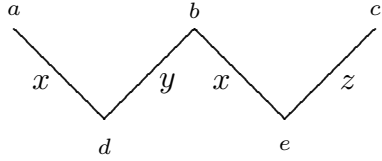
{\beginpicture
    \setcoordinatesystem units <1cm,1cm>
  \multiput{} at  0 0  4 1 /
  \plot 0 1  1 0  2 1  3 0  4 1 /
  \multiput{$x$} at 0.3 0.4  2.3 0.4 /
  \multiput{$y$} at 1.7 0.4 /
  \multiput{$z$} at 3.7 0.4 /
  \put{$\ssize a$} at 0 1.2
  \put{$\ssize b$} at 2 1.2
  \put{$\ssize c$} at 4 1.2
  \put{$\ssize  d$} at 1 -.3
  \put{$\ssize  e$} at 3 -.3

\endpicture} 
$$
\vglue-.2cm
$\s$
	\medskip
{\bf 5.3. Type 3: At least three pairwise non-isomorphic bristles.}
	\medskip
Proof. Let $B_1,B_2,B_3$ be pairwise non-isomorphic subbristles of the special $K(3)$-module $W$.
According to 3.2, we have $B_i\cap B_j = 0$ for $i\neq j.$
Let $x', y, z'$ be non-zero elements of $E$ such that $x'(B_2+B_3) =  y(B_1+B_3) = z'(B_1+B_3) = 0.$
Since $B_1$ and $B_3$ are non-isomorphic, $B_1+B_3 = B_1\oplus B_3$ 
contains only the bristles $B_1,B_3,$ thus
$B_1+B_2+B_3 = W$. If we would have 
$y\in \langle x',z'\rangle,$ then $y$ would annihilate $W = B_1+B_2+B_3$, but
we assume that $W$ faithful. Thus $x',y,z'$ is a basis of $E$.
Let $a\in B_1\setminus \soc B_1,\ b\in B_2\setminus \soc B_2,\ 
c\in B_3\setminus \soc B_3.$ Since $B_1+B_3 = B_1\oplus B_3,$ the elements $x'a,\ z'c$ form a basis
of $W_1.$ Therefore $yb = \alpha x'a + \gamma z'c$ for some (non-zero) elements $\alpha,\gamma\in k.$
Let $x = \alpha x'$ and $z = \gamma z',$ so that $yb = xa+zc.$ 
We get the following coefficient quiver:
$$
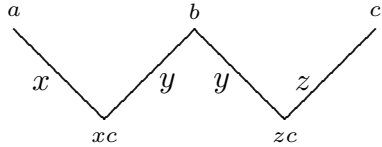
{\beginpicture
    \setcoordinatesystem units <1cm,1cm>
  \multiput{} at  0 0  4 1 /
  \plot 0 1  1 0  2 1  3 0  4 1 /
  \multiput{$x$} at 0.3 0.4 /
  \multiput{$y$} at 1.7 0.4  2.3 0.4 /
  \multiput{$z$} at 3.2 0.4 /
  \put{$\ssize a$} at 0 1.2
  \put{$\ssize b$} at 2 1.2
  \put{$\ssize c$} at 4 1.2
  \put{$\ssize xc$} at 1 -.2
  \put{$\ssize zc$} at 3 -.2
\endpicture} 
$$
\vglue-.2cm
$\s$

	\medskip
{\bf 5.4. Type $\infty.$ At least two isomorphic bristles.}
	\medskip 
Proof. 
Let $N$ be a bar which is the direct sum of two isomorphic bristles, say bristles
isomorphic to $B$. Let $B = (k,E/\langle y,z\rangle),$ where $x,y,z$ is a basis of $E$.
Then $\Ker x$ is not contained in $N$, since otherwise we get an embedding of $S_0$ into $W$.
Let $c\in \Ker x \setminus N$. Then $c$ does not generate a bristle (since otherwise $W$
would not be faithful), thus $yc, zc$ is a basis of $W_1.$ There are elements $a,b\in N_0$
with $xa = zc$ and $xb = yc$. Altogether we get
$$
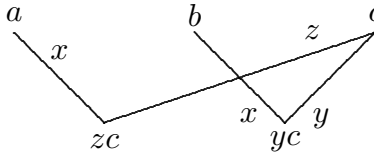
{\beginpicture
    \setcoordinatesystem units <1cm,1cm>
  \multiput{} at  0 0  4 1 /
  \plot 0 1  1 0  4 1  3 0  2 1  /
  \multiput{$x$\strut} at 0.5 0.8  2.6 0.1 /
  \multiput{$z$} at  3.3 1 /
  \multiput{$y$\strut} at  3.4 0.1 /

  \put{$a$} at 0 1.2
  \put{$b$} at 2 1.2
  \put{$c$} at 4 1.2
  \put{$zc$} at 1 -.2
  \put{$yc$} at 3 -.2
\endpicture} 
$$
\vglue-.2cm
$\s$

	\medskip
{\bf 5.5. Remark 1.} We have seen: In case the bristle type of $W$ is not $1$, then $W$ is a
tree module. Let us add: {\it A special $K(3)$-module with
precisely one bristle is not a tree module.}
	\medskip
Proof. There are only two possibilities for the coefficient quiver of a tree module 
with dimension vector $(3,2)$,
namely
$$
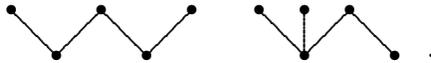
{\beginpicture
    \setcoordinatesystem units <.5cm,.5cm>
\put{\beginpicture
\multiput{$\ssize \bullet$} at 0 1  1 0  2 1  3 0  4 1 /
\plot 0 1  1 0  2 1  3 0  4 1 /
\endpicture} at 0 0
\put{\beginpicture
\multiput{$\ssize \bullet$} at 0 1  1 0  1 1  2 1  3 0 /
\plot 0 1  1 0  2 1  3 0 /
\plot 1 1  1 0 /
\endpicture} at 5 0 
\put{.} at 7.3 -.5
\endpicture}
$$
In both cases, we see that there are at least 2 bristles.
(Actually, only the first case can occur, since in the second case, $W$ contains a submodule
with dimension vector $(2,1)$, and such a submodule is either uniform or has $S_0$ as a 
submodule.)
$\s$ 
	\medskip
{\bf Remark 2.} All the coefficient quivers which occur are without parameters. This is
clear for tree modules. In the remaining case (type 1), 
there is just one bar, and this bar
contains a bristle. Thus the submodule $U$ given by the cycle is not a bar. The fact that
$U$ is faithful implies immediately that any parameter can be deleted.
	\bigskip\bigskip

{\bf 6. $L(3)$-modules.}
	\medskip
We want to adapt Propositions 3.7 and 3.8, as well as 2.5  to $L(3)$-modules. 
As we have mentioned in the introduction, $L(3)$
is the 4-dimensional local $k$-algebra with radical square zero, 
thus $L(3) = k\langle x,y,z\rangle/(x,y,z)^2.$ The radical of $L(3)$ will be denoted by $E$.
The dimension vector $\bdim M$ of an $L(3)$-module $M$
is the pair of numbers $(\dim M/EM,\dim EM).$
	\medskip
{\bf 6.1. Atoms and bristles.} The notation to be used corresponds to the notation for
$K(3)$-modules introduced in section 2.1;
the precise correspondence will be stated in section 6.3.
	\medskip 
If $0 \neq a\in E$, then $C(a) = C(\langle a \rangle) = {}_AA/\langle a \rangle$ is 
an atom, any atom is obtained in this way, and 
{\it $a \mapsto C(a)$ provides a bijection between
$\Bbb PE$ and the set of isomorphism classes of atoms.}
If $a_1,a_2$ are linearly independent in $E$, then $B(a_1,a_2) = B(\langle a_1,a_2\rangle) =
{}_AA/\langle a_1,a_2\rangle$
is a bristle, any bristle if obtained in this way, and 
{\it $\langle a_1,a_2\rangle$ provides a
bijection between the set of lines in $\Bbb PE$ and the set of isomorphism classes of bristle.}
If $0\neq a \in E$, and $l$ is a line in $\Bbb E$, then $\langle a \rangle \in l$ iff
$C(a)$ maps onto $B(l).$
	\medskip
The relation between $K(3)$ and $L(3)$ is given by the
push-down functor $\pi\:\mod K(3) \to \mod L(3)$ which will be described in 
two different (of course related) ways (partly, we follow
the presentation in the appendix A of [RZ3]).  
	\medskip
{\bf 6.2. The push-down functor $\pi$.} This is the functor
$$
 \pi\:\mod K(3) \to \mod L(3)
$$ 
which sends 
the $K(3)$-module $M = (M_0,M_1,\phi\:E\otimes M_0 \to M_1)$ 
to the $L(3)$-module $\pi M = (M_0\oplus M_1,\left[\smallmatrix 0 &\phi\cr 0 & 0 
\endsmallmatrix\right])$. 
The decisive property of $\pi$ is the following: If $M,M'$ are 3-Kronecker modules without
non-zero simple direct summands, then 
we have
$$
 \Hom_{L(3)}(\pi(M),\pi(M')) = 
 \pi\Hom_{K(3)}(M,M') \oplus
 \Hom_{L(3)}(\pi(M),\pi(M'))_{ss};
$$
here, we denote by 
$\Hom_{L(3)}(X,X')_{ss}$ the set of $L(3)$-homomorphisms
with semisimple image, for arbitrary $L(3)$-modules  $X, X'$. 
	\medskip
{\it The functor $\pi$ provides a bijection between the isomorphism classes of the
indecomposable $K(3)$-modules $M$ different from $S_1$ and the isomorphism classes
of the indecomposable $L(3)$-modules, such that $\bdim M = \bdim \pi M$.} (Of course, also
$\pi S_1$ is defined; but $\bdim S_1 = (0,1)$, whereas $\bdim \pi S_1 = (1,0$.)
	\medskip
The $L(3)$-module ${}_{L(3)}L(3)$ is just $\pi(k,E)$. Of course, the radical of $\pi(k,E)$ is
$0\oplus E = E.$ Let us stress that
in this way, {\it the arrow space $E$ of the quiver $K(3)$ is identified with
the radical of $L(3).$} 
	\medskip
There is an {\bf alternative way} to define $\pi$.  
The $K(3)$-modules may be seen as the $R$-modules, where $R$ is the path algebra of the
quiver $K(3)$, thus $R = 
\left[\smallmatrix k & 0 \cr E & k\endsmallmatrix\right]$ is the ring of all
matrices of the form $\left[\smallmatrix \lambda & 0 \cr e & \lambda'\endsmallmatrix\right]$,
with $\lambda, \lambda'\in k$ and $e\in E$. Here, the $K(3)$-module 
$(M_0,M_1;\phi_M)$
corresponds to the set $\left[\smallmatrix M_0 \cr M_1\endsmallmatrix\right]$ of 
column matrices, and the module multiplication is just 
matrix multiplication using the map $\phi_M$.
By abuse of notation, we denote the ring $R$ also by $K(3)$ (this explains our
convention to call the $K(3)$-modules just $K(3)$-modules).
We may consider $L(3)$ as the subring of $K(3)$ given by the matrices of the form 
$\left[\smallmatrix \lambda & 0 \cr e & \lambda\endsmallmatrix\right]$ with 
$\lambda\in k,$ and $e\in E$. 
Then, the functor $\pi$ is the restriction functor, which sends the 
$K(3)$-module $M$ to the $L(3)$-module often denoted by $M|_{L(3)}$ (the
map $K(3) \to \End_k M$ given by the module multiplication is restricted to $L(3)$). 
In particular, in this interpretation, 
the underlying vector spaces of $M$ and $\pi M$ are the same.
This interpretation is particular helpful when we compare the submodule structure of
$M$ and of $\pi M,$ see 6.3 and 6.4.
	\medskip
{\bf 6.3. Socle properties.} {\it Let $M$ be a $K(3)$-module. Always, we 
have $|\soc M| = |\soc \pi M|.$
If $M$ is indecomposable and not simple, then $\soc M = M_1 = J(\pi M) = \soc \pi M.$}
	\medskip
Proof. It is enough to look at indecomposable $K(3)$-modules. If $M$ is simple, then also
$\pi M$ is simple. Now assume that $M = (M_0,M_1,\phi_M)$ is indecomposable and not simple.
Since $M$ is not isomorphic to $S_0$, we have
$\soc M = M_1$. 
On the other hand, $\pi M$ is an indecomposable $L(3)$-module of Loewy length
$2$, therefore $\soc \pi M = J(\pi M)$ and $J(\pi M)$ is the image of $\phi_M$, thus
$J(\pi M) = M_1.$ $\s$
	\medskip
{\bf Corollary.} {\it The $K(3)$-module $M$ is uniform iff the $L(3)$-module $\pi M$ is uniform.}
	\medskip
{\bf Atoms.} 
{\it The functor $\pi$ provides a bijection between the isomorphism classes
of the atom $K(3)$-modules  and the isomorphism classes of the atom $L(3)$-modules with
$$
 \pi \widetilde C(a) = C(a)
$$
where $0 \neq a \in E$.} 
	\medskip 
{\bf Bristols.} 
{\it The functor $\pi$ provides a bijection between the isomorphism classes
of the bristle $K(3)$-modules  and the isomorphism classes of the bristle $L(3)$-modules with
$$
 \pi \widetilde B(a_1,a_2) = B(a_1,a_2)
$$
where $a_1,a_2$ are linearly independent in $E$.}

	\medskip
{\bf 6.4. Submodules.}
Given a $K(3)$-module $M$,
its $K(3)$-submodules are, of course, also $L(3)$-submodules, but usually, there are 
$L(3)$-submodules of $\pi M$ which are not $K(3)$-submodules. 
There is the following general observation:
	\medskip
{\bf Lemma.} 
{\it Let $M$ be a $K(3)$-module. Then any $L(3)$-submodule $N$ of $\pi M$ with $M_1 \subseteq N$
is a $K(3)$-submodule of $M.$}
	\medskip
Proof. We note that $K(3)$ is generated as a vector space by $L(3)$ and the matrix 
$p = \left[\smallmatrix 0 & 0 \cr 0 & 1 \endsmallmatrix\right]$, and we have 
$pM = M_1.$ Thus 
any $L(3)$-submodule $N$ of $M$ with $M_1\subseteq N$ is a $K(3)$-submodule.
$\s$
	\medskip
{\bf 6.5. Special $K(3)$-modules, special $L(3)$-modules.}
	\medskip
An $L(3)$-module $M$ will be said to be {\it special} provided $\bdim M = (3,2),$
$M$ is faithful, and has no simple direct summand and no uniform submodule of length 3.
	\medskip
{\bf Lemma.} {\it Let $W$ be a $K(3)$-module with $\bdim W = (3,2).$ Then $W$ is a special 
$K(3)$-module iff $\pi W$ is a special $L(3)$-module.}
	\medskip
Proof. 
According to condition (ii) in 2.3, $W$ is a special $K(3)$-module iff $W$ is
faithful, indecomposable and has no uniform submodule of length 3.

Let $M$ be an arbitrary $K(3)$-module. As we have mentioned, $M$ is 
indecomposable iff $\pi M$ is indecomposable. Also, $M = (M_0,M_1)$ is faithful iff there is no
non-zero element $e\in E$ with $eM_0 = \phi_M(e\otimes M) = 0$, iff there is no non-zero
element $e\in E$ with $e(\pi M) = 0$ iff $\pi M$ is faithful. Also, if $U$ is a uniform
submodule of $M$, then $U$ is an $L(3)$-submodule of $\pi M$, and $U$ is uniform also 
as an $L(3)$-module.

Conversely, assume that $V$ is a uniform $L(3)$-submodule of $\pi M$. Let $V = \pi U$.
The embedding $v\:V \to \pi M$ 
is of the form $v = \pi u + v'$, where $u\:U \to M$ is $K(3)$-linear
and $v':U \to $ is an $L(3)$-homomorphism with semisimple image. 
Then $v'\:V \to \pi M$ 
vanishes on $\soc V$. Thus $v = \pi u + v'$ shows that $u$
does not vanish on $\soc U$. Since $U$ is uniform, $u$ is a monomorphism, thus $M$ has
a uniform submodule of length 3.
$\s$
		\medskip
Following [RZ3], an $L(3)$-module $M$ is said to be {\it solid}
provided $M = \pi W$ for some $K(3)$-module $W$ with $\End W = k.$
	\medskip
{\bf Corollary.} {\it An $L(3)$-module $M$ is special iff $M$ has dimension vector
$(3,2)$, is faithful and solid.}
	\medskip
We are interested in bristle and bar submodules of these special $L(3)$-modules. 
	\medskip
{\bf 6.6. Bristle submodules.}
{\it Let $M$ be a $K(3)$-module. 
The bristle $K(3)$-submodules are bristle $L(3)$-submodules of $\pi M$ and for
any bristle $L(3)$-submodule $U'$ of $\pi M$, there is a bristle $K(3)$-submodule $U$ of $M$
such that $\soc U' = \soc U.$
All bristle $L(3)$-submodules $U'$ and $U''$ of $\pi M$ 
with $\soc U' = \soc U''$ are isomorphic.}
	\medskip
The following formulation provides more details: Let $R = K(3)$ and $A = L(3) \subset K(3)$.
Let $B = Rb$ be a subbristle of ${}_RM$, then 
$B$ is a subbristle of $\pi M$ with $B\cap M_0 \neq 0$. 
For any element $x\in M_1$, $A(b+x)$ is a subbristle of $\pi M$,
and the $R$-submodule generated by $A(b+x)$ is $Rb+Rx$ (thus a bristle iff $Rx \subseteq Rb$).
{\it We have $A(b+x) = A(b+x')$ iff $x-x' \in B\cap M_1;$ 
also, $A(b+x) \cap M_0 \neq 0$ iff $x\in B\cap M_1.$}
Conversely, any subbristle of $\pi M$ is of the form $A(b+x)$, where
$Rb$ is an $R$-bristle and $x\in M_1.$
	\medskip
Bristle modules are defined for any algebra, bar modules not. Similar to the case of a
$K(3)$-module with dimension vector $(3,2)$, we consider here an $L(3)$-module $M$ with
dimension vector $(3,2)$ and call a submodule $N$ of $M$ a {\it bar}, 
provided $\bdim N = (2,2)$ and $N$ is not faithful.
	\medskip 
{\bf 6.7. Bar submodules of special modules.} {\it Let $W$ be a special $K(3)$-module. 
Any bar $L(3)$-submodule of $\pi W$ is actually a $K(3)$-submodule.} 
	\medskip
Proof. This follows directly from 6.4.
$\s$
	\medskip

{\bf 6.8. Proposition.} 
{\it A special $L(3)$-module has at most $3$ isomorphism classes of 
subbristles and at most $3$ bar submodules. Any subbristle of a special $L(3)$-module
is contained in a bar submodule.}
	\medskip
{\bf Remark.} One has to be aware that there is a decisive difference between special 
$K(3)$-modules and special $L(3)$-modules: Whereas a special $K(3)$-module may have only finitely
many subbristles (and then 1, 2, or 3), a special $L(3)$-modules always has infinitely
many subbristles, see 3.9 and 6.6. 
	\medskip
Proof of Proposition 6.8.
Let $W$ be a special $K(3)$-module. The bar $K(3)$-submodules are
the bar $L(3)$-submodules of $\pi W.$ 
The bristle $K(3)$-submodules are bristle $L(3)$-submodules of $\pi W$ and for
any bristle $L(3)$-submodule $U$ of $\pi W$, there is a bristle $K(3)$-submodule $U'$ of $W$
such that $U$ and $U'$ are isomorphic $K(3)$-modules and $\soc U = \soc U'.$
Altogether we see that 6.8 is a direct consequence of 3.7 and 3.8.
$\s$
	\medskip
{\bf 6.9.} If $M$ is a special $L(3)$ module, the number of simple modules which occur as
the socle of a subbristle of $M$ will be called the {\it bristle type} of $M$.
{\it If $M$ is a special $K(3)$-module, then the bristle 
type of $M$ is equal to the bristle type of $\pi M.$}
	\medskip
{\bf 6.10.} Finally, we show that also the analog of 2.5 holds true.
Note that the assertion (ii$'$) has to be modified, see 6.11.
	\medskip
{\bf Proposition.}
{\it Let $M$ be a special $L(3)$-module. Let $0 \neq a \in E.$ The following
conditions are equivalent:

\item{\rm(i)} The atom $C(a)$ has no factor module which is a subbristle of $M$.
\item{\rm(ii)} The trace of $C(a)$ in $M$ is an atom.
\item{\rm(ii$'$)} The trace of $C(a)$ in $M$ is indecomposable and has length $3$. \par}
	\medskip
Proof. Let $M = \pi W$, where $W$ is a special $K(3)$-module. 

(i) $\implies$ (ii). We assume that $C(a)$ has no factor module which is a subbristle of $M$. 
If $\phi\:\widetilde C(a) \to W$ would be a $K(3)$-homomorphism whose image is a bristle, then
$\pi(\phi)$ would be an $L(3)$-homomorphism whose image is a bristle. 
Thus, we see that 
$\widetilde C(a)$ has no factor module which is a subbristle of $W$, this is condition (i)
in 2.5. It follows that the trace $D$ of $\widetilde C(a)$ in $W$ is an atom. Of course,
$\pi(D)$ is a submodule of $M$ which is an atom.
First of all, this means that there is a homomorphism $\psi\:\widetilde C(a) \to W$
with image $D$. As a consequence, $\pi(\psi)\:C(a) \to M$ has image in $\pi D$. 
Second, any homomorphism $\widetilde C(a) \to W$ maps into $D$. 
Let $f\:C(a) \to M$ be a homomorphism. According to 6.2, we can write $f = \pi(\psi')+f'$,
where $\psi':\widetilde C(a) \to W$ and where $f'\:C(a) \to M$ has semisimple image. 
As we have mentioned, $\psi'$ maps into $\pi(D)$. Since $\pi(D)$ is an atom, its socle has
length 2, and therefore coincides with the socle of $M$. Since the image of $f'$ is semisimple,
this image is contained in $\soc M = \soc \pi D$, thus in $\pi D$. This shows that
$\pi D$ is the trace of $C(a)$ in $M$. 

(ii$'$) $\implies$ (i). We assume that the trace $N$ 
of $C(a)$ in $M$ has length 3 and is indecomposable.
Since $M$ has no uniform submodule of length 3, it follows that $N$ is an atom. 
There is $f\:C(a) \to M$ whose image is $N$. 
According to 6.3, we write $f = \pi(\phi)+f'$, with a $K(3)$-homomorphism 
$\phi\:\widetilde C(a) \to  W$ and an $L(3)$-homomorphism $f':C(a) \to M$ with semisimple
image. The image of $f'$ lies in $\soc M$. Since the image of $f$ does not lie in $\soc M$,
the image of $\phi$ cannot lie in $\soc W$. Therefore the image of $\phi$ has length 2 or 3.
Now, the image of $\phi$ cannot be a bristle, since otherwise also the image of $\pi(\phi)$ would
be a bristle. This shows that the image of $\phi$ has length 3. Let $D$ be the trace of
$\widetilde C(a)$ in $W$. Then $\Im(\phi) \subseteq D$, thus $|D| \ge 3.$ 
On the other hand, the trace $D$ of $\widetilde C(a)$ in $W$ is mapped under $\pi$
into the trace $N$ of $C(a)$ in $M$. But $\pi(D) \subseteq N$ implies that $|D| \le |N| = 3.$
This shows that the trace $D$ of $\widetilde C(a)$ in $W$ has length 3. 
According to 2.5, it follows that 
$\widetilde C(a)$ has no factor module which is a subbristle of $W.$ 

Assume that $C(a)$ has a factor module $B$ which is a subbristle of $M$. 
Since there is an $L(3)$-epimorphism $C(a) \to B$, there is also a $K(3)$-epimorphism 
$\widetilde C(a) \to \widetilde B$, where $\pi(\widetilde B) = B.$
Since $B$ is a subbristle of $M$, $\widetilde B$ is a subbristle of $W$, thus
$\widetilde C(a)$ has a factor module which is a subbristle of $W,$ a contradiction.
$\s$
	\medskip
{\bf 6.11. Remark.} {\it If $M$ is a special $L(3)$-module and $0\neq a\in E,$
the trace of $C(a)$ in $M$ has length $3$ or $4.$ In both cases, the trace of $C(a)$ 
in $M$ may be decomposable.}
	\medskip
Proof. Let $W$ be a special $L(3)$-module and $M = \pi W.$
Since $L(3)$ is a local algebra and $C(a)$ is non-zero, the 
socle $\soc M$ of $M$ is contained in the trace $D$ of $C(a)$ in $M$. Also, there is a non-zero
homomorphism $\phi\:\widetilde C(a) \to W.$ Since 
the image of $\pi(\phi)\:C(a) \to M$ is not contained in $\soc M$, 
we see that $\soc M$ is a proper submodule of $D$. Thus, $|D| \ge 3.$ 
On the other hand, $D$ is not faithful, thus $D$ is a proper submodule of $M$, 
therefore $|D| \le 4.$ 

Here is an example where the trace has length 3 and
is decomposable: Let $W$ be the special $K(3)$-module
of bristle type 1, as exhibited in Proposition 4.2. The module $\widetilde C(y)$
is the local module annihilated by $y$, it has as factor module the bristle $B = B(y,z)$
annihilated by $y$ and $z$. Now $B$ is the submodule of $W$ generated by $(100)$, 
and this is the trace of $\widetilde C(y)$
in $W$. The trace $D$ of $C(y)$ in $M = \pi W$ is $\pi B + \soc M$. Thus $D$ has length 3
and is the direct sum of $\pi B$ and a simple module.

Here is an example where the trace has length 4 and
is decomposable: Let $W$ be the special $K(3)$-module
of bristle type 2, as exhibited in Proposition 4.2. The trace of the module $\widetilde C(y)$
in $W$ is the direct sum of two bristles, thus also the trace of $C(y)$
in $M = \pi W$ is the direct sum of two bristles.
$\s$
	\bigskip\bigskip
{\bf 7. The special algebras.}
	\medskip
{\bf 7.1. Lemma.} {\it Let $A$ be a short local algebra with radical $J$ and assume that
$J^2 \neq 0.$ The
following conditions are equivalent.
\item{\rm (i)} ${}_AJ$ is faithful as an $A/J^2$-module. 
\item{\rm (i$'$)} If $z\in J$ with $zJ = 0,$ then $z\in J^2.$
\item{\rm (i$''$)} $\soc J_A \subseteq J^2.$ 
\item{\rm (ii)} $J^2 = \soc J_A$. 
\item{\rm (iii)} $J_A$ has no simple direct summand.\par}
	\medskip
Proof. The equivalence of (i) and (i$'$) is just the definition of faithfulness.
(i$'$) implies (i$''$): The socle of a right module $M_A$ are just the elements $m\in M$
with $mJ = 0.$ 

(i$''$) implies (ii):
Since $A$ is short, we have $J^2 \subseteq \soc J_A$.

(ii) implies (iii):
Since $J^2 \ne 0,$ we know that $J_A$ is a right module of Loewy length $2$.
A right module $M_A$ of Loewy length 2 has no simple direct summand iff $MJ = \soc M_A$. 

(iii) implies (i$'$). Assume there is $z\in J\setminus J^2$ with $zJ = 0.$ Then
$\langle z\rangle$ is a simple direct summand of $J_A.$
$\s$ 
	\medskip
{\bf 7.2. Lemma.} {\it Assume that $A$ is a short local algebra of Hilbert type $(3,2)$ with
$J^2 = \soc {}_AA = \soc A_A.$ Then the following conditions are
equivalent:
\item{\rm (i)} ${}_AA$ has no uniform submodule of length $3$.
\item{\rm (ii)} $A$ has no serial module of length $3$.
\item{\rm (i$^*$)} $A_A$ has no uniform submodule of length $3$.
\item{\rm (ii$^*$)} $A^{\op}$ has no serial module of length $3$.\par}
	\medskip
Proof. (i) $\implies$ (ii): Assume that there exists a serial module $M$ of length $3$. 
Since $M$ is local, its projective cover is ${}_AA$ and $\Omega M$
is a left ideal of length 3. The socle of $\Omega M$ is contained in $\soc {}_AA$, thus of
length 1 or 2. If the socle of $\Omega M$ would be of length 2, 
then $\soc M = \soc {}_AA =
J^2,$ and then $M = {}AA/\Omega M$ has Loewy length at most 2, a contradiction. 

(ii) $\implies$ (i): Assume that ${}_AA$ has a uniform submodule $U$ of length 3. Then
$J^2 = \soc {}_AA$ is not contained in $U$, therefore ${}_AA/U$ is a module of length 3 which
is of Loewy length 3. Of course, a module of length 3 and Loewy length 3 is serial.

The equivalence of (ii) and (ii$^*$) is given by $k$-duality.
$\s$ 
	\medskip
{\bf 7.3. Proposition.} {\it Let $A$ be a short local algebra of
  Hilbert type $(3,2).$ The following conditions are equivalent.
\item{\rm (i)} 
   $J^2 = \soc {}_AA = \soc A_A$ 
   and there is no uniform left ideal with length $3.$
\item{\rm (i$'$)} 
   $J^2 = \soc {}_AA = \soc A_A$ 
   and there is no uniform right ideal with length $3.$
\item{\rm (ii)} 
   ${}_AJ$ is a special $L(3)$-module.
\item{\rm (ii$'$)} 
   $J_A$ is a special $L(3)$-module.
\item{\rm (iii)} ${}_AJ$ is a solid $A$-module, $J_A$ is a solid
   $A^{\op}$-module. \par}
	\medskip
Recall that we call an algebra {\it special} provided it satisfies the equivalent conditions
mentioned in Proposition 7.3.
It is shown in [RZ3] that for a short local algebra with 
Hilbert type $(3,2)$ the existence of a non-projective reflexive module
implies that both ${}_AJ$ and $J_A$ are solid, thus that $A$ is special.
This is the implication (i) $\implies$ (iii) in the main theorem
of the present paper. 

The importance of condition (ii) lies in the fact that this condition 
refers only to properties of $A$ as a left $A$-module. 
	\medskip
Proof of Proposition 7.3. The equivalence of (i) and (i$'$) is shown in 7.2. 

The equivalence of (i) and (ii) is seen as follows: To say that ${}_AJ$ is special means that
${}_AJ$ is faithful as an $A/J^2$-module, that ${}_AJ$ has no simple direct summand and
no uniform submodule of length 3. According to 7.1, the condition $J^2 = \soc J_A$ is equivalent to
the condition that ${}_AJ$ is faithful as an $A/J^2$-module. Second, $J^2 = \soc {}_AJ$ is
equivalent to the condition that ${}_AJ$ has no simple direct summand. 

Since the conditions (i) and (i$'$) are equivalent, also (ii) and (ii$'$) are equivalent.
According to Corollary 6.5, (ii) implies that ${}_AJ$ is solid. Since 
(ii) and (ii$'$) are equivalent, we see that these conditions imply that also $J_A$ is solid.
Thus (ii) implies (iii).

Finally, let us show that (iii) implies (ii). We
assume that both ${}_AJ$ and $J_A$ are solid. 
Since $J_A$ is solid, it has no simple direct summand, thus, according to 7.1, ${}_AJ$ is faithful
as an $A/J^2$-module (that means, as an $L(3)$-module).  Since ${}_AJ$ is also solid, Corollary 6.5 asserts that
${}_AJ$ is special.  
$\s$
	\medskip
Let us mention an important property of these special algebras.
	\medskip 
{\bf 7.4. Lemma.} {\it Let $A$ be a special algebra. Then any
left or right ideal of $A$ of length $3$ or $4$ is a two-sided ideal.}
	\medskip
Proof. Since the opposite of a special algebra is a special algebra, we only have to look at left
ideals.

Since $A$ is short, $J^2 \subseteq \soc{}_AJ \cap \soc J_A.$ 
Since ${}_AJ$ is solid, it has no simple direct summand, thus $\soc {}_AJ \subseteq J^2.$
Since ${}_AJ$ is faithful as an $A/J^2$-module, 7.1 asserts that $J_A$ has no simple
direct summand, thus also $\soc J_A \subseteq J^2.$ Altogether, we see that 
$J^2 = \soc{}_AA = \soc A_A.$

Since $A$ is special, we know that ${}_AJ$ has no uniform submodule of length 3. 
It follows that any left ideal of $A$ of length at least 3 contains $\soc{}_AJ = J^2.$
Thus, let $U$ be a left ideal of $A$ included in $J$. We want to see that $Ua \subseteq U$
for any $a\in A.$ Now $a = \lambda+a'$ for some $\lambda \in k$ and $a'\in J.$
Of course, $U\lambda \subseteq U$. And $Ua' \subseteq J^2 \subseteq U.$
$\s$
	\medskip
Since for a special algebra $A$ the radical $J$ considered as a left module is a special
$L(3)$-module, and any proper left ideal of $A$ is contained in $J$, we can apply Proposition 6.8
and obtain:
	\medskip

{\bf 7.5. Proposition.} {\it A special algebras has at most $3$ isomorphism classes of bristle modules
which occur as left ideals and at most $3$ left ideals which are bars. Any bristle left ideal
is contained in a left ideal which is a bar.}  $\s$
	\medskip
{\bf 7.6.} If $A$ is a special algebra with radical $J$ such that 
${}_AJ =\pi W$, then we call the 
bristle type of $W$
the {\it left bristle type} of $A$; it
is the number of simple left ideals of $A$ which occur as socle of a bristle left ideal.
Similarly, the {\it right bristle type of} $A$ is the number of 
simple right ideals of $A$ which occur as socle of a bristle right ideal.
	\bigskip\bigskip
{\bf 8. The correspondence between bristle left ideals and right bars.}
	\medskip
{\bf 8.1. Lemma.} {\it Let $A$ be a special algebra. 
If $B$ is a left module which is a bristle, let $l(B) = \{a\in A\mid aB = 0\}$, this is the
left annihilator of $B$. Then $\Omega(_AB) = l(B)$, and this is a $4$-dimensional ideal of $A$ and $B = A/l(B).$ In this way, we obtain a
bijection between the isomorphism classes of bristles and the 
$4$-dimensional ideals of $A$.}
	\medskip
Proof. Let $M$ be a local module and 
$f\:{}_AA \to M$ be a projective cover, thus $\Omega({}_AM) = \Ker f.$
Of course, $l(M) = \{a\in A\mid aM = 0\} \subseteq \Ker f = \Omega({}_AM),$ since 
for $a\in l(M)$ we have $f(a) = f(a\cdot 1) = af(1) = 0.$
Actually, $l(M)$ is the largest twosided ideal of $A$ which is contained in $\Ker f$
(namely, if $I$ is a twosided ideal of $A$, then $I(A/I) = 0$). 

Now, if $M = B$ is a left bristle, then $\Omega({}_AB)  $ is 4-dimensional 
and any 4-dimensional left ideal of $A$ is a twosided ideal, see 7.4. This shows
that $\Omega(_AB) = l(B)$. It follows that $B \simeq {}_AA/l(B).$ 
Thus, if $B,B'$ are bristles with $l(B) = l(B'),$ then $B$ and $B'$ are isomorphic.

Of course, conversely, isomorphic modules have the same annihilator. 
Thus, if $B, B'$ are isomorphic left modules, then 
$l(B) = l(B').$ 
$\s$
	\medskip
{\bf 8.2. Lemma.} {\it Let $A$ be a special algebra. 
Under this bijection the isomorphism classes of the bristles which occur as left ideals 
correspond bijectively to the right bars.}
	\medskip
Proof. Let $B = Aa$ be a left bristle.
Then $\Omega({}_AB) = l(B)$ is annihilated from the right by $a$, thus
not faithful as right $A/J^2$-module. Conversely, assume that $N_A$ is a bar of $J_A$, say 
annihilated by $a\in J\setminus J^2$ (from the right). 
Consider the map ${}_AA \to {}_AA$ given by right multiplication 
with $a$; its image is $Aa$, its kernel is $\Omega({}_AAa).$ Since $a\in J\setminus J^2$, 
$Aa$ has length at least 2.
Since $Na = 0$, we see that $N$ is contained in the kernel. But $N$ is of dimension 4, thus
$Aa$ has length 2, and $N = \Omega({}_AAa).$
	\medskip
{\bf Remark.} The bijection is between {\bf isomorphism classes} of bristles and
{\bf individual} ideals. In particular, this has to be noticed in case we consider left ideals
which are bristles and right bars, thus on the one hand we deal with 
left ideals of dimension 2 (better: with the corresponding isomorphism classes) and on the other
hand with ideals of dimension 4.
	\medskip
{\bf 8.3.} 
Applying Lemma 8.1 to $A^{\op}$, we similarly get {\it a bijection between the isomorphism classes of
right bristles and the $4$-dimensional ideals of $A$,} where a right bristle $B$ is sent to
$r(B) = \{s\in A\mid Ba = 0\}$, whereas a $4$-dimensional ideal $I$ of $A$ is sent to $B = A/I$,
considered as a right $B$-module. 

But note: since $L(3)$ is commutative, any left bristle $B$ may be considered also as a right
bristle and conversely; also, the left annihilator $l(B)$ coincides with the right annihilator $r(B)$
and one obtains in this way just the $4$-dimensional ideals. 
	\medskip

{\bf 8.2$^*$. Lemma.} {\it Let $A$ be a special algebra. 
Under the bijection between the bristle
 right modules and the $4$-dimensional ideals of $A$
the isomorphism classes of the bristles which occur as right ideals correspond bijectively to the 
left bars.} $\s$

	\bigskip\bigskip
{\bf 9. Proof of Theorem 1.1.}
	\medskip
We assume that $A$ is a special algebra.  
We fix a complement $E$ of $J^2$ in $J$, thus $E$ is a subspace of $J$ with $J = J^2\oplus E.$
Of course, we may identify $E$ with $J/J^2$ (thus with the radical of $L(3) = A/J^2$).
	\medskip
{\bf 9.1. The projective plane $\Bbb PE.$}
{\it The projective plane $\Bbb PE$ as index set for the  atoms.}
If $a\in J\setminus J^2,$ let  $U(a) = Aa+J^2$, thus we have $C(a) = {}_AA/U(a).$
Similarly, {\it The set of lines in $\Bbb PE$ is index set for the bristles.}
If $\langle a\rangle $ is an element of $\Bbb PE$, and $l$ is a line in $\Bbb PE$, then 
$\langle a\rangle\in l$
iff $B(l)$ is a factor module of $C(a).$ 

As we have seen in Lemma 7.4, any left or right ideal of
$A$ of length 3 or 4 is a twosided ideal. Thus, the left bristles may be considered also as right bristles, the left atoms as right atoms. 
	\medskip
{\bf 9.2. The lines $g(B)$ and $t(N).$} 
If $N$ is a 4-dimensional ideal of $A$, let $t(N) = \{\langle a \rangle\mid a\in N\cap E\}.$
Clearly, {\it $t(N)$ is a line in $\Bbb PE$} (it is the line joining $a$ and $a'$, where
$a,a'$ are elements of $E$ such that $N = Aa+Aa'$). 
	\medskip 
If $B$ is a bristle, we write $g(B)$ for the set of all elements
$\langle a \rangle\in \Bbb PE$, where $a\in E$,  such that $B$ is a factor module of $C(a).$
Of course, $g(B)$ is a line in $\Bbb PE$, and all lines in $\Bbb PE$ are
obtained in this way, 
see 6.1. 
	\medskip
{\bf Lemma.} {\it Let $A$ be a special algebra. If $B$ is a bristle, then }
$$
 g(B) = t(\Omega_AB).
$$
	\medskip
Proof. First, consider $a\in E$ with $\langle a \rangle\in g(B).$ 
Thus there is an epimorphism $\phi\:C(a) \to B.$ Since $C(a)$ is local,
the projective cover of $C(a)$ is ${}_AA$ and
there is a commutative diagram with exact rows
$$
\CD
 0 @>>> Aa+\soc A @>>> {}_AA @>>> C(a) @>>> 0 \cr
 @.     @V\mu VV            @|             @VV\phi V \cr
 0 @>>> \Omega_AB @>>> {}_AA @>>> B @>>> 0 .
\endCD
$$
with an inclusion map $\mu.$ The inclusion map $\mu$ shows that $a\in \Omega_AB \cap E$
is a non-zero element and thus $\langle a \rangle\in t(\Omega_AB).$

Conversely, if $0\neq a\in \Omega_AB \cap E$ (so that $\langle a \rangle \in t(\Omega_AB)$),  
then $Aa+\soc A \subseteq \Omega_AB$,
and we get the left part of the commutative diagram seen above (with $\mu$ the inclusion map).
Completing it on the right,
we obtain the required epimorphism $C(a) \to B.$
$\s$
	\medskip
{\bf 9.3. Proposition.} {\it Let $A$ be a special algebra. 
Let $0\neq a\in E$. The following conditions are equivalent:
\item{\rm (i)} The element $\langle a \rangle$ in $\Bbb PE$ 
  does not belong to a left bar, and $aA$ is not a right bristle. 
\item{\rm (ii)} The atom $C(a)$ is extensionless (and $\Omega C(a)$ is an atom).\par}
	\medskip
Remark. The bracket in (ii) adds a general fact: 
If $M$ is an indecomposable non-projective module which is extensionless, then
$\Omega M$ is indecomposable, again, see [RZ3]. 
Of course, $C(a)$ is indecomposable and non-projective. Thus, if
$C(a)$ is extensionless, then $\Omega C(a)$ is indecomposable. Since $C(a) = {}_AA/U(a),$ 
with $U(a) = Aa+J^2$, we have $\Omega C(a) = U(a)$. If $U(a)$ is indecomposable, then it is an atom.
	\medskip
Proof of 9.3. (i) implies (ii). First, let us assume that $a$ does not belong 
to a left bar. Then $Aa$ is not a bristle (see 7.5),
thus an atom. Let $D$ be the trace of $U(a)$ in ${}_AJ,$ thus $U(a)\subseteq D.$
Since $D$ is generated by $Aa$ and $Aa$ is not faithful as an $L(3)$-module,
also $D$ is not faithful.
Since ${}_AJ$ is faithful as an $L(3)$-module, 
$D$ is a proper submodule of ${}_AJ$.
If $D$ has length 4,then $D$ is a left bar. But $a$ is not contained in a left bar, 
$D$ has length 3 and therefore $D = U(a)$.
Since $aA$ is not a bristle, it follows that $aA = U(a)$. As a consequence, the 
embedding $Aa = U(a) \to {}_AJ$ is a left ${}_AA$-approximation, of course a minimal one.
This shows that $0 \to Aa \to {}_AA \to C(a) \to 0$ is an $\mho$-sequence, thus $C(a)$
is extensionless. 

(ii) implies (i). Assume that $C(a)$ is extensionless. Then $\Omega C(a) = U(a)$ is 
indecomposable, thus equal to $Aa$, and $0 \to Aa \to {}_AA \to C(a) \to 0$ is an $\mho$-sequence.

Assume that $a$ belongs to a left bar $N$. Since $Aa$ is a local module of length 3, and $N$ is
not faithful as an $L(3)$-module, 
$Aa$ generates the module $N$, thus there is a homomorphism $\phi\:Aa \to N$
with image not contained in $Aa.$ Since the inclusion map $\mu\:Aa \to {}_AA$ is a left
${}_AA$-approximation, there is $\phi'\:{}_AA \to {}_AA$ with $\phi = \phi'\mu.$
Now $\phi'$ is the right multiplication by some element $b\in A$, therefore 
$\phi(a) = \phi'\mu(a) = ab\in aA.$ Thus $\phi(Aa) = A\phi(a) \subseteq AaA = U(a) = Aa.$
This contradiction shows that $a$ does not belong to a left bar. 

Finally, we show that $U(a) = aA$. Let $x \in \soc U(a).$ The map $\phi\:Aa \to {}_AA$ 
with $\phi(a) = x$ factors through the 
inclusion map $\mu\:Aa \to {}_AA$, thus
there is $\phi'\:{}_AA \to {}_AA$ with $\phi = \phi'\mu.$
Again, $\phi'$ is the right multiplication by some element $b\in A$, therefore 
$x = \phi(a) = \phi'\mu(a) = ab\in aA.$ This shows that $\soc U(a) \subseteq aA,$ therefore
$aA = U(a).$ In particular, $aA$ is not a bristle.
$\s$
	\medskip
{\bf 9.4. Proposition.} {\it Let $A$ be a special algebra. 
Let $0\neq a\in E$. The following conditions are equivalent:
\item{\rm (i)} The element $\langle a \rangle$ in $\Bbb PE$ does not belong to a right bar.
\item{\rm (i$'$)} The atom $C(a)$ has no factor module which is a left subbristle.\par}
\item{\rm (ii)} {\it The trace of $C(a)$ in ${}_AA$ is an atom.}
	\medskip
Let us add: {\it If the trace of $C(a)$ in ${}_AA$ is an atom, then $C(a)$ is torsionless.}
The converse is of course not true. 
	\medskip
Proof of Proposition. The equivalence of (i) and (i$'$) follows from Lemma 9.2. which asserts
that $g(B) = t(\Omega B)$.
For the equivalence of (i$'$) and (ii), see 6.10 with $M = {}_AJ.$ 
$\s$
	\medskip
{\bf 9.5. Corollary.} {\it If $a$ does not belong to a left bar nor to a right bar, then
$C(a)$ is both extensionless and torsionless.} Thus, $\Omega C(a)$ is a reflexive atom.
	\medskip
Proof. Assume that $a$ does not belong to a left bar nor to a right bar. 
Then $aA$ is not a right bristle, since otherwise $a$ belongs to a right bar, see 7.5.
According to Proposition 9.3, $C(a)$ is extensionless 
and $\Omega C(a) = U(a)$ is an atom.
According to Proposition 9.4, $C(a)$ is torsionless. 
Also, $\mho U(a) = \mho\Omega C(a)$, thus $\mho U(a) = C(a)$.
This shows that $\mho U(a)$ is torsionless. Since both $U(a)$ and $\mho U(a)$
are torsionless, $U(a)$ is reflexive.
$\s$
	\medskip
{\bf 9.6. Proof of Theorem 1.1.} We take an element $0\neq a\in E$ which does not belong to
a left bar nor to a right bar. Such elements do exist, since there are only finitely many
left bars and finitely many right bars, and $k$ is algebraically closed.

Since $\langle a \rangle$ does not belong to a right bar, $aA$ is not a
right bristle. According to 9.5, $\Omega C(a)$ is a reflexive atom.
$\s$
	\medskip
{\bf 9.7. Remark.} {\it Let $A$ be a short local algebra with radical $J$.
Let $U$ be an atom submodule of $J$.
Even if the trace of $U$ in ${}_AJ$ is an atom,
the embedding $U \to J$ does not have to be a left ${}_AA$-approximation.}
	\medskip
{\bf Example.}
Let $A$ be the algebra shown in the upper row on the left (with $A^{\op}$ on the right). 
Let $U = Ay.$
In the lower row, we depict  $\mho\Omega C(y) = \mho (Ay)$.
$$
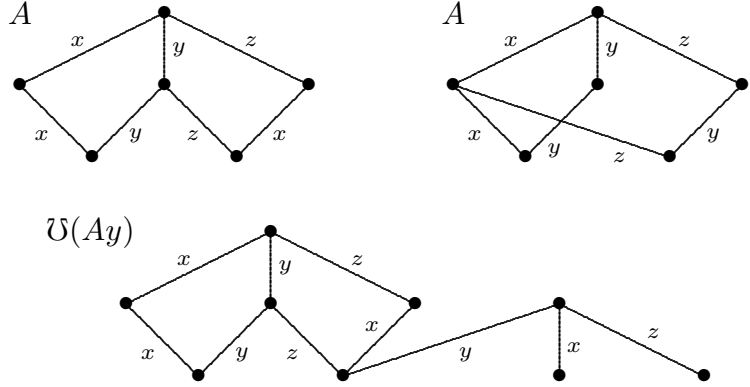
{\beginpicture
    \setcoordinatesystem units <.8cm,.8cm>
\put{\beginpicture
\put{$A$} at 0 2
\multiput{$\bullet$} at 0 1  2 1  4 1  1 0  3 0  2 2 /
\plot 2 2  0 1  1 0  2 1  3 0  4 1  2 2 /
\plot 2 2  2 1 /
\multiput{$\ssize x$} at 0.8 1.6  0.3 0.3  3.6 0.3 /
\multiput{$\ssize y$} at 2.2 1.5  1.6 0.3 /
\multiput{$\ssize z$} at 3.2 1.6  2.4 0.3 /
\endpicture} at 0 0
\put{\beginpicture
\put{$A$} at 0 2
\multiput{$\bullet$} at 0 1  2 1  4 1  1 0  3 0  2 2 /
\plot 2 2  0 1  1 0  2 1 /
\plot 0 1  3 0  4 1  2 2 /
\plot 2 2  2 1 /
\multiput{$\ssize x$} at 0.8 1.6  0.3 0.3 /
\multiput{$\ssize y$} at 2.2 1.5  1.4 0.1  3.6 0.3 /
\multiput{$\ssize z$} at 3.2 1.6  2.3 0 /
\endpicture} at 6 0
\put{\beginpicture
\put{$\mho (Ay)$} at -.5 2
\multiput{$\bullet$} at 0 1  2 1  4 1  1 0  3 0  2 2  6 1  6 0  8 0 /
\plot 2 2  0 1  1 0  2 1  3 0  4 1  2 2 /
\plot 2 2  2 1 /
\plot 3 0  6 1  6 0 /
\plot 6 1  8 0 /
\multiput{$\ssize x$} at 0.8 1.6  0.3 0.3  3.4 0.7  6.2 0.4 /
\multiput{$\ssize y$} at 2.2 1.5  1.6 0.3  4.7 0.3 /
\multiput{$\ssize z$} at 3.2 1.6  2.3 0.3  7.3 0.6 /
\endpicture} at 3 -3 
\endpicture}
$$
Note that for all pairs $(\beta,\gamma)\neq 0$, the right ideal 
$(\beta y+\gamma z)A$ is a bristle, thus $C(\beta y+\gamma z)$
is not extensionless. ({\it The atom $C(\alpha x+\beta y+\gamma z)$ is 
extensionless iff $\alpha\neq 0$ and $\gamma \neq 0$}.) The picture shows the case
$\alpha= 0,\ \beta = 1,\ \gamma = 0.$
	\medskip
{\bf 9.8.} As we have seen in Proposition 9.3 and Remark 9.7
the fact that $C(a  ) = A/Aa$ is not extensionless
depends not only on the left module structure of the radical $J$, but also on the 
right module structure of $J$. 

Given a special algebra $A$ with radical $J = \pi W$, where $W$ is a (special)
$K(3)$-module, and $U$ an atom in $W$, 
there is a 
short local algebra $A'$ with the ``same'' radical $\pi W$ (where ``same'' refers
to the left-module structure) such that $\pi U$, considered as an $A'$-module, 
is Gorenstein-projective. Namely:
	\medskip
{\bf Proposition.} 
{\it Let $\bdim W = (3,2).$ Let $U$ be a submodule of $W$ with $\bdim U = (1,2)$
such that also the trace of $U$ in $W$ is an atom. 

Assume that the annihilator of $U$ in $E$ is generated by $x\in E.$ 
There is a short local algebras $A$
with an element $x\in J$ such that ${}_AJ = W$ and $U = Ax.$ 

If $A$ is a
short local algebra with an element $x\in J$ such that ${}_AJ = W$ and $U = Ax,$ then
$x$ is a (left and right) Conca element of $A$} (thus, $U$ is a Gorenstein projective
$A$-module with $\Omega U \simeq U$).
	\medskip 
Recall that if $A$ is a short local algebra, then $x\in A$ is said to be a (left and right) 
{\it Conca element}
provided $x^2 = 0$ and $xJ = Jx = J^2$ (see [AIS] and also [RZ3]).
	\medskip
Proof. Let $a,b,c$ be a generating set of $W$ such that $x$ generates $U$. Let $A$ be
generated by $x,y,z,$ with radical $J$ such that $x\cdot1 = a,\ y\cdot 1 = b,\  z\cdot 1 = c.$
$\s$
	\bigskip\bigskip 

{\bf 10. The left layout and the right layout of a special algebra.}
	\medskip
{\bf 10.1. Lemma.
} {\it The right layout of a special algebra is isomorphic to the left layout.}
	\medskip
Proof. 
The following table shows in the left column the various special algebras $A$ 
and the corresponding  opposite algebra $A^{\op}$ in the right column. We assume that 
$A$ is a special algebra with radical $\pi W$, where $W$ is a
$K(3)$-module with coefficient quiver as displayed in section 4. Note that this means
that we work with a fixed basis $x,y,z$ of $E$, and also with a fixed basis (say $(100), (010),
(001)$) of $W_0$. In our identification of ${}_AJ$ with $\pi W$, the elements $(100), (010), (001)$
of $\pi W$ will correspond (in some order) to a second basis $x',y',z'$ of $E$, the order
which we choose is seen in the pictures in the left column. 

It is straightforward to observe that indeed 
the opposite algebra has the presentation as shown
in the right column. In particular, its radical (and this is the right module $J_A$) has a
coefficient quiver which is obtained from the coefficient quiver we started with by replacing
$x,y,z$ by $x',y',z'$.

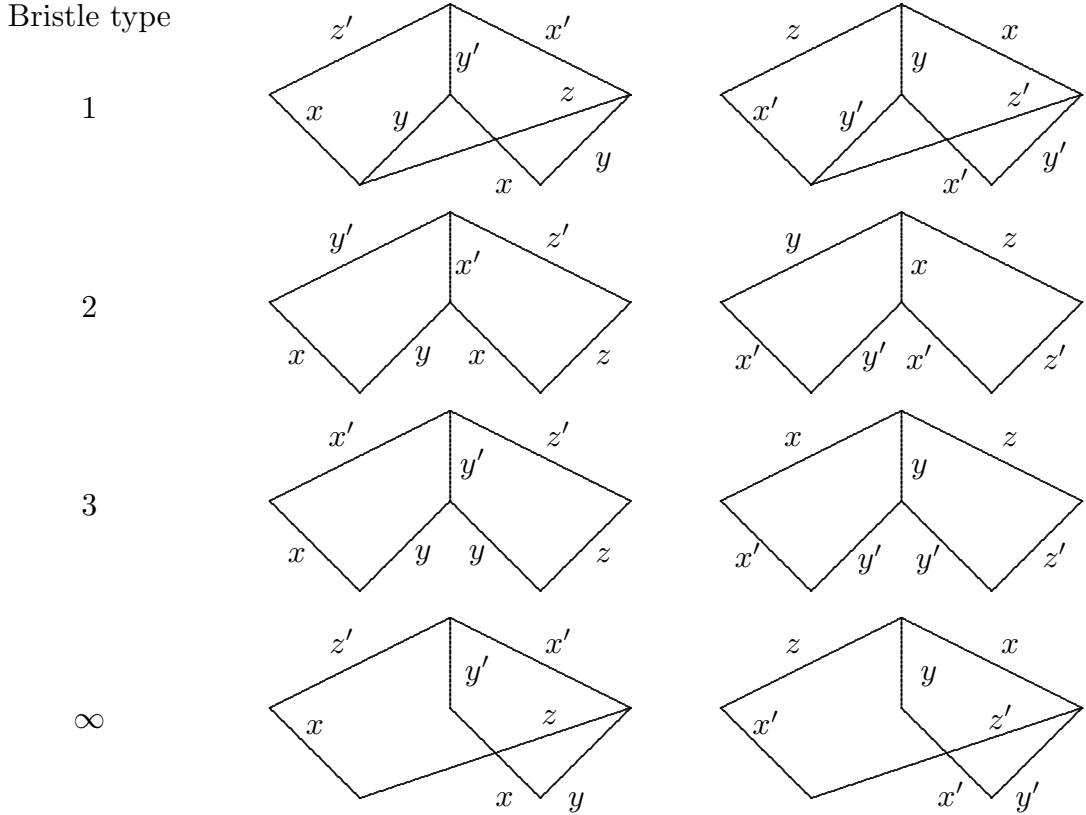
$$
{\beginpicture
    \setcoordinatesystem units <1cm,1cm>
\put{Bristle type} at 3 5
\put{\beginpicture
  \multiput{} at  0 0  4 1 /
  \plot 0 1  1 0  2 1  3 0  4 1  1 0 /
  \multiput{$x$\strut} at 0.5 0.8  2.6 0 /
  \multiput{$z$} at  3.3 1 /
  \multiput{$y$\strut} at 1.45 0.75  3.7 0.3 /

  \plot 0 1  2 2  4 1 /
  \plot 2 2  2 1 /
  \put{$z'$\strut} at  0.8 1.7
  \put{$y'$\strut} at  2.2 1.4
  \put{$x'$\strut} at  3.2 1.7
\endpicture} at 7 4
\put{\beginpicture
  \multiput{} at  0 0  4 1 /
  \plot 0 1  1 0  2 1  3 0  4 1  1 0 /
  \multiput{$x'$\strut} at 0.5 0.8  2.6 0 /
  \multiput{$z'$} at  3.3 1 /
  \multiput{$y'$\strut} at 1.45 0.75  3.7 0.3 /

  \plot 0 1  2 2  4 1 /
  \plot 2 2  2 1 /
  \put{$z$\strut} at  0.8 1.7
  \put{$y$\strut} at  2.2 1.4
  \put{$x$\strut} at  3.2 1.7
\endpicture} at 12 4
\put{1} at 3 4

\put{\beginpicture
  \multiput{} at  0 0  4 1 /
  \plot 0 1  1 0  2 1  3 0  4 1 /
  \multiput{$x$} at 0.3 0.4  2.3 0.4 /
  \multiput{$y$} at 1.7 0.4 /
  \multiput{$z$} at 3.7 0.4 /

  \plot 0 1  2 2  4 1 /
  \plot 2 2  2 1 /
  \put{$y'$\strut} at  0.8 1.7
  \put{$x'$\strut} at  2.2 1.4
  \put{$z'$\strut} at  3.2 1.7
\endpicture} at 7 1.8
\put{\beginpicture
  \multiput{} at  0 0  4 1 /
  \plot 0 1  1 0  2 1  3 0  4 1 /
  \multiput{$x'$} at 0.3 0.4  2.2 0.4 /
  \multiput{$y'$} at 1.7 0.4 /
  \multiput{$z'$} at 3.7 0.4 /

  \plot 0 1  2 2  4 1 /
  \plot 2 2  2 1 /
  \put{$y$\strut} at  0.8 1.7
  \put{$x$\strut} at  2.2 1.4
  \put{$z$\strut} at  3.2 1.7
\endpicture} at 12 1.8
\put{2} at 3 1.8

\put{\beginpicture
  \multiput{} at  0 0  4 1 /
  \plot 0 1  1 0  2 1  3 0  4 1 /
  \multiput{$x$} at 0.3 0.4 /
  \multiput{$y$} at 1.7 0.4  2.3 0.4 /
  \multiput{$z$} at 3.7 0.4 /

  \plot 0 1  2 2  4 1 /
  \plot 2 2  2 1 /
  \put{$x'$\strut} at  0.8 1.7
  \put{$y'$\strut} at  2.25 1.4
  \put{$z'$\strut} at  3.2 1.7
\endpicture} at 7 -.4
\put{\beginpicture
  \multiput{} at  0 0  4 1 /
  \plot 0 1  1 0  2 1  3 0  4 1 /
  \multiput{$x'$} at 0.3 0.4 /
  \multiput{$y'$} at 1.65 0.35  2.3 0.35 /
  \multiput{$z'$} at 3.7 0.4 /

  \plot 0 1  2 2  4 1 /
  \plot 2 2  2 1 /
  \put{$x$\strut} at  0.8 1.7
  \put{$y$\strut} at  2.2 1.4
  \put{$z$\strut} at  3.2 1.7
\endpicture} at 12 -.4
\put{3} at 3 -.4

\put{\beginpicture
  \multiput{} at  0 0  4 1 /
  \plot 0 1  1 0  4 1  3 0  2 1  /
  \multiput{$x$\strut} at 0.5 0.8  2.6 0 /
  \multiput{$z$} at  3.1 .9 /
  \multiput{$y$\strut} at  3.4 0 /
  \plot 0 1  2 2  4 1 /
  \plot 2 2  2 1 /
  \put{$z'$\strut} at  0.8 1.7
  \put{$y'$\strut} at  2.3 1.4
  \put{$x'$\strut} at  3.2 1.7

\endpicture} at 7 -2.8
\put{\beginpicture
  \multiput{} at  0 0  4 1 /
  \plot 0 1  1 0  4 1  3 0  2 1  /
  \multiput{$x'$\strut} at 0.5 0.8  2.55 0 /
  \multiput{$z'$} at  3.1 .9 /
  \multiput{$y'$\strut} at  3.4 0 /
  \plot 0 1  2 2  4 1 /
  \plot 2 2  2 1 /
  \put{$z$\strut} at  0.8 1.7
  \put{$y$\strut} at  2.3 1.4
  \put{$x$\strut} at  3.2 1.7

\endpicture} at 12 -2.8
\put{$\infty$} at 3 -2.8
\endpicture}
$$
$\s$
	\medskip
{\bf 10.2.}
As we have seen, the left layout and the right layout have the same shape.
Of course,
as subsets of $\Bbb PE$ they may be different. Let us present an example with
3 left subbristles and 3 right subbristles, all pairwise different.
	\medskip
{\bf Example.} {\it  A special algebra $A$ 
with six isomorphism classes of bristles which are left or right subbristles.}
The algebra $A$ is defined by the following relations: 
$$
  zx = 0,\ x^2 = yx,\ xy = y^2 = zy,\ xz = 0,\ yz = z^2,\ y^2+x^2+z^2.
$$
$$
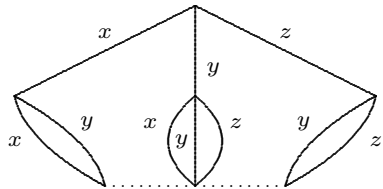
{\beginpicture
    \setcoordinatesystem units <1cm,1cm>
\multiput{} at -2 0  2 2 /
\plot -2 1  0 2  2 1 /
\plot 0 2  0 1 /
\setquadratic
\plot -2 1 -1.3 0.5  -1 0 /
\plot -2 1 -1.7 0.5  -1 0 /

\plot 0 1  -0.3 0.5  0 0 /
\plot 0 1   0   0.5  0 0 /
\plot 0 1   0.3 0.5  0 0 /

\plot 2 1 1.3 0.5  1 0 /
\plot 2 1 1.7 0.5  1 0 /
\setlinear
\setdots <1mm>
\plot -1 0  1 0 /

\multiput{$\ssize x$} at -1 1.7   -2 0.5  -.5 0.7 /
\multiput{$\ssize y$} at 0.2 1.3  -1.2 0.7  1.2 0.7  -.15 0.5 /
\multiput{$\ssize z$} at 1 1.7   2 0.5  .45 0.7 /
\endpicture}
$$
One easily checks that
$A$ is a special algebra. Also:

\item{$\bullet$} $Ax$ is annihilated by $z$ and $x-y$; 
\item{$\bullet$} $Ay$ is annihilated by $x-y$ and $y-z$;
\item{$\bullet$} $Az$ is annihilated by $x$ and $y-z.$
\item{$\bullet$} $(x-y)A$ is annihilated by $x$ and $y$;
\item{$\bullet$} $(y-z)A$ is annihilated by $y$ and $z$;
\item{$\bullet$} $(x-y+z)A$ is annihilated by $x$ and $z$.
	\smallskip
Therefore: {\it 
The left ideals 
$$
  Ax,\ Ay,\ Az
$$ 
are bristles, and the right ideals 
$$
 (x-y)A,\ (y-z)A,\ (x-y+z)A
$$ 
are also bristles,} thus there are 6 isomorphism classes of bristles which occur
as left or right ideals. Also, {\it there are 6 four-dimensional ideals which are left bars 
or right bars, namely}
$$
\gather
 Ax+Ay,\ Ax+Az,\ Ay+Az;\cr
  (x-y)A+(y-z)A,\ (x-y)A+zA,\ (y-z)A+xA.
\endgather
$$
		\bigskip
{\bf 11. Commutative algebras.}
	\medskip
{\bf 11.1. Proposition.} {\it Let $M$ be a special $L(3)$-module. There is a commutative
algebra $A$ such that ${}_AJ = M$ as $A/J^2$-module.}
	\medskip
Proof. Here are such algebras (on the right, we write down the defining relations, in order
to obtain $A$ from $k[x,y,z]$):
$$
{\beginpicture
    \setcoordinatesystem units <1cm,1cm>
\put{Bristle type} at 3 5
\put{Relations} at 12 5
\put{\beginpicture
  \multiput{} at  0 0  4 1 /
  \plot 0 1  1 0  2 1  3 0  4 1  1 0 /
  \multiput{$x$\strut} at 0.5 0.8  2.6 0.1 /
  \multiput{$z$} at  3.3 1 /
  \multiput{$y$\strut} at 1.45 0.75  3.7 0.3 /

  \plot 0 1  2 2  4 1 /
  \plot 2 2  2 1 /
  \put{$z$\strut} at  0.8 1.7
  \put{$y$\strut} at  2.2 1.4
  \put{$x$\strut} at  3.2 1.7
\endpicture} at 7 4
\put{1} at 3 4
\put{$yz,\ x^2,\ y^2-xz,\ z^2$} at 12 4
\put{\beginpicture
  \multiput{} at  0 0  4 1 /
  \plot 0 1  1 0  2 1  3 0  4 1 /
  \multiput{$x$} at 0.3 0.4  2.3 0.4 /
  \multiput{$y$} at 1.7 0.4 /
  \multiput{$z$} at 3.7 0.4 /

  \plot 0 1  2 2  4 1 /
  \plot 2 2  2 1 /
  \put{$y$\strut} at  0.8 1.7
  \put{$x$\strut} at  2.2 1.4
  \put{$z$\strut} at  3.2 1.7
\endpicture} at 7 1.8
\put{2} at 3 1.8
\put{$xz,\ yz,\ y^2,\ x^2-z^2$} at 12 1.6

\put{\beginpicture
  \multiput{} at  0 0  4 1 /
  \plot 0 1  1 0  2 1  3 0  4 1 /
  \multiput{$x$} at 0.3 0.4 /
  \multiput{$y$} at 1.7 0.4  2.3 0.4 /
  \multiput{$z$} at 3.7 0.4 /

  \plot 0 1  2 2  4 1 /
  \plot 2 2  2 1 /
  \put{$x$\strut} at  0.8 1.7
  \put{$y$\strut} at  2.2 1.4
  \put{$z$\strut} at  3.2 1.7
\endpicture} at 7 -.4
\put{3} at 3 -.4
\put{$xy,\ xz,\ yz,\ x^2 \! +\! y^2\! +\! z^2$} at 12 -.6

\put{\beginpicture
  \multiput{} at  0 0  4 1 /
  \plot 0 1  1 0  4 1  3 0  2 1  /
 \multiput{$x$\strut} at 0.5 0.8  2.6 0 /
  \multiput{$z$} at  3.1 .9 /
  \multiput{$y$\strut} at  3.4 0 /

  \plot 0 1  2 2  4 1 /
  \plot 2 2  2 1 /
  \put{$z$\strut} at  0.8 1.7
  \put{$y$\strut} at  2.2 1.4
  \put{$x$\strut} at  3.2 1.7

\endpicture} at 7 -2.8
\put{$\infty$} at 3 -2.8
\put{$xy,\ x^2,\ y^2,\ z^2$} at 12 -2.8
\endpicture}
$$
	\medskip
When dealing with a commutative special algebra, a bar is a 4-dimensional ideal which is
not faithful as an $L(3)$-module.  
	\medskip
{\bf 11.2. Proposition.} {\it Let $A$ be a commutative special algebra and $a\in E$.
The following conditions are equivalent:
\item{\rm (i)} The element $\langle a \rangle$ does not belong to a bar.
\item{\rm (ii)} The atom $C(a)$ is extensionless.
\item{\rm (iii)} The trace of $C(a)$ in ${}_AA$ is an atom (in particular, $C(a)$ is torsionless).
\item{\rm (iv)} The atom $C(a)$ is Gorenstein-projective. 
\item{\rm (v)} The atom $C(a)$ is reflexive.
\par}
	\smallskip
{\it If these conditions are satisfied, then $C(a)$ is $\Omega$-periodic with period at most $2$.}
	\medskip
Concerning the condition (iii): Assume that the trace of $C(a)$ is an atom $Y$.
Since $Y$ is local, there is a homomorphism $\phi\:C(a) \to {}_AA$ with image $Y$. Since 
both $C(a)$ and $Y$ have length 3, $\phi$ is an isomorphism. 
This shows that $C(a)$ is torsionless. Remark 11.6 will show that the converse is not true:
$C(a)$ may be torsionless, whereas its trace in ${}_AA$ is not an atom.
	\medskip
Proof of 11.2. The equivalence of (i) and (ii) follows from Proposition 9.3, the equivalence
of (i) and (iii) follows from Proposition 9.4.

Let us show that the equivalent conditions (i), (ii), (iii)
imply (iv). Thus, assume that $a$ does not belong to a bar. According to (ii),
$C(a)$ is extensionless, and therefore $U(a) = \Omega C(a)$ is indecomposable, thus
equal to $Aa$. According to (iii), the trace of $C(a)$ in ${}_AA$ is an atom, say equal
to $Ab$ for some $b\in E.$ 

For any $x\in A$, let $\rho_x\:{}_AA \to {}_AA$ be the 
right multiplication by $x$. Then $Ab = U(b)$ is the image of $\rho_b$ and 
$Aa$ is the image of $\rho_a.$ Look at the composition
$$
    {}_AA @>\rho_a>> {}_AA @>\rho_b>> {}_AA @,
$$
The image of $\rho_a$ is $Aa = U(a)$. The image $Ab$ of $\rho_b$ is the 
trace of $C(a)$ in ${}_AA$. It follows that $U(a)$ is the kernel of $\rho_b.$
In particular, we have $ba = 0$ and the sequence exhibited above is exact. 

Since $A$ is commutative, we also have $ab = 0.$ Therefore, also the sequence
$$
    {}_AA @>\rho_b>> {}_AA @>\rho_a>> {}_AA 
$$
is a complex. Now. the image $Ab$ of $\rho_b$ has length 3. Also the image $Aa$ of
$\rho_a$ has length 3, thus the kernel of $\rho_a$ has length 3, and therefore this
sequence is also exact. Altogether, we see that $U(b) = Ab = \Omega Aa$. It follows
that $C(a) = Ab$ is $\Omega$-periodic with period at most 2. This proves the 
additional assertion. 

If we consider the complex $P_\bullet$
$$
  \cdots @>>> P_2 @>\rho_b>> P_1 @>\rho_a>> P_0 @>\rho_b>> P_{-1} @>>> \cdots
$$
with all modules $P_i = {}_AA$ and 
using alternatively $\rho_a$ and $\rho_b$, we obtain an exact complex. 
The $A$-dual of $P_\bullet$ is again $P_\bullet$, since the $A$-dual of $\rho_x\:{}_AA \to {}_AA$
is just the left multiplication by $x$.
Thus shows that $C(a)$ is the image in an exact complex $P_\bullet$ whose $A$-dual is also
exact, thus $C(a)$ is Gorenstein-projective.

Of course, (iv) implies (ii), and also (iv) implies (v).
$\s$
	\medskip
(v) implies (iii).
Assume that $C(a)$ is reflexive. Then $C(a)$ and $\mho C(a)$ both are
torsionless and we have an $\mho$-sequence $0 \to C(a) \to {}_AA^t \to \mho C(a) \to 0,$
for some $t\ge 1.$ 
If $t \ge 2$, then the Loewy length of the cokernel 
of $C(a) \to {}_AA^t$ would be 3, but this cokernel is $\mho C(a)$ and $\mho C(a)$ has
Loewy length 2, since it is torsionless (and has no non-zero projective direct summand).
Therefore $t = 1.$ Since the embedding $C(a) \to {}_AA$ is a left ${}_AA$-approximation, 
the image of $C(a)$ in ${}_AA$ has to be the trace of $C(a)$ in ${}_AA$, thus the trace of
$C(a)$ in ${}_AA$ is an atom. 
$\s$
	\medskip
{\bf 11.3. Remark.} In condition (iii) of 11.2, it is not enough to assume that 
the atom $C(a)$ is torsionless:
{\it The module $C(a)$ may be torsionless, without being Gorenstein projective.}
As an example, take the following commutative algebra $A$.
	
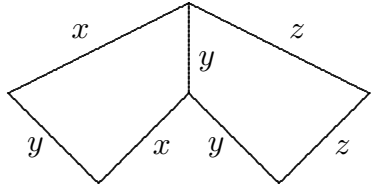
\medskip
$$
{\beginpicture
    \setcoordinatesystem units <1cm,1cm>
\put{\beginpicture
  \multiput{} at  0 -.2  4 1 /
  \plot 0 1  1 0  2 1  3 0  4 1 /
  \multiput{$y$} at 0.3 0.4  2.3 0.4 /
  \multiput{$x$} at 1.7 0.4 /
  \multiput{$z$} at 3.7 0.4 /

  \plot 0 1  2 2  4 1 /
  \plot 2 2  2 1 /
  \put{$x$\strut} at  0.8 1.7
  \put{$y$\strut} at  2.2 1.4
  \put{$z$\strut} at  3.2 1.7
\endpicture} at 7 -.4
\endpicture}
$$
Here, the atoms $C(x)$ and $C(z)$ both are torsionless, but $\Omega C(x)$ and $\Omega C(z)$
are decomposable. Thus, for $a = x$ and $a = z$, the module $C(a)$ is torsionless, 
but not extensionless, and, as a consequence, not Gorenstein-projective.
The traces of $C(x)$ and $C(z)$ in ${}_AA$ have length 4 (the trace of $C(x)$
is decomposable, the trace of $C(z)$ is indecomposable).
$\s$
	\medskip
{\bf 11.4. Remark.} In case $W$ is a special $K(3)$-module of bristle
type $\infty$, there is also
an anti-commutative short local algebra $A$ with ${}_AJ = \pi W$ as $A/J^2$-module:
$$
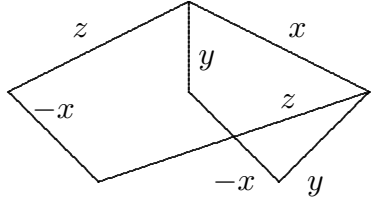
{\beginpicture
    \setcoordinatesystem units <1cm,1cm>
\put{\beginpicture
  \multiput{} at  0 0  4 1 /
  \plot 0 1  1 0  4 1  3 0  2 1  /
 \multiput{$-x$\strut} at 0.5 0.8  2.5 0 /
  \multiput{$z$} at  3.1 .9 /
  \multiput{$y$\strut} at  3.4 0 /

  \plot 0 1  2 2  4 1 /
  \plot 2 2  2 1 /
  \put{$z$\strut} at  0.8 1.7
  \put{$y$\strut} at  2.2 1.4
  \put{$x$\strut} at  3.2 1.7

 \endpicture} at 7 -2.8

\endpicture}
$$
	\bigskip\medskip
{\bf 12. More about $6$-dimensional local algebras.}
	\medskip
Let $A$ be a local algebra which is not self-injective and has a non-projective
reflexive module. We have mentioned in the introduction that in case $A$ is 
short the dimension of $A$ has to be at least 6. Using Proposition 4.2, one may
construct many isomorphism classes of short 
local algebras of dimension 6 which are not self-injective, but have non-projective
reflexive modules. Let us stress that there do exist also 6-dimensional local algebras
which are not self-injective, with non-projective reflexive modules, and which are not short.
	\medskip
{\bf 12.1.} {\it 
$6$-dimensional local algebras with non-projective reflexive modules need not to
be self-injective or short.}
	\medskip
Example: Let $A' = k\langle x,y\rangle/(x^2,y^2,yxy)$ and $M = A'x = A'/A'x.$
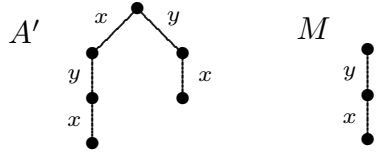
$$
{\beginpicture
    \setcoordinatesystem units <.5cm,.5cm>
\put{\beginpicture
\put{$A'$} at -1.5 2.5
\multiput{$\bullet$} at 0 0  0 1  0 2  1 3  2 1  2 2 /
\plot 0 0  0 2  1 3  2 2  2 1 /
\multiput{$\ssize x$} at 0.2 2.8  -.4 0.5  2.5 1.5 /
\multiput{$\ssize y$} at  -.4 1.5  1.8 2.8 /
\endpicture} at 0 0 
\put{\beginpicture
\put{} at 0 3
\put{$M$} at -1.2 2.5
\multiput{$\bullet$} at 0 0  0 1  0 2  /
\plot 0 0  0 2  /
\multiput{$\ssize x$} at  -.4 0.5  /
\multiput{$\ssize y$} at  -.4 1.5  /
\endpicture} at 5 0
\endpicture}
$$
{\it The module $M$ is $\Omega$-periodic with period $1$, and Gorenstein projective,
thus reflexive.}
	
	\bigskip\medskip
{\bf Acknowledgment.} 
This paper is a continuation of our joint work with Zhang Pu from SJTU Shanghai devoted to 
the study of modules over short local algebras. Unfortunately, due to the severe Corona restrictions,
the cooperation had to be suspended. The investigations presented here strongly rely on the
previous joint papers and the author has to acknowledge his dependence on Pu's insight into the realm of Gorenstein-projective modules. In addition, the 
author wants to thank two referees for spotting misprints and giving helpful advice.
	\bigskip\medskip

{\bf 13. References.}
\frenchspacing
	\medskip
\item{[AIS]} L. L. Avramov, S. B. Iyengar, L. M. \c Sega. 
    Free resolutions over short local rings. J.
    London Math. Soc. 78 (2008), 459-–476.
\item{[Ra]} M. Ramras.
    Betti numbers and reflexive modules, Ring theory (Proc. Conf., Park City, Utah,
    1971), Academic Press, New York, 1972, 297–-308.
\item{[R1]} C. M. Ringel. Representations of k-species and bimodules. J. Algebra 41 (1976),
    269–-302.
\item{[R2]} C. M. Ringel. Exceptional modules are tree modules.
    Lin. Alg. Appl. 275-276 (1998), 471--493. 
\item{[R3]} C. M. Ringel. The elementary 3-Kronecker modules. arXiv:1612.09141.
\item{[RZ1]} C. M. Ringel, P. Zhang.
    Gorenstein-projective and semi-Gorenstein-projective modules. 
    Algebra \& Number Theory 14--1 (2020), 
    1–36. dx.doi.org/10.2140/ant.2020.14.1
\item{[RZ2]} C. M. Ringel, P. Zhang.
     Gorenstein-projective and semi-Gorenstein-projective modules. II.
    Journal of Pure and Applied Algebra. Vol 224, Issue 6 (2020).
    Article Number: 106248, https://doi.org/10.1016/j.jpaa.2019.106248 
\item{[RZ3]} C. M. Ringel, P. Zhang.  
    Gorenstein-projective modules over short local algebras.
    Journal of the London Mathematical Society, Vol 106, Issue 2 (2022), 528--589.
    \newline
    http://dx.doi.org/10.1112/jlms.12577

	\bigskip\smallskip
{\baselineskip=1pt
\rmk
C. M. Ringel\par
Fakult\"at f\"ur Mathematik, Universit\"at Bielefeld \par
POBox 100131, D-33501 Bielefeld, Germany  \par
ringel\@math.uni-bielefeld.de\par
\medskip
ORCID: 0000-0001-7434-0436
}

	\bigskip\bigskip
{\bf Statements and Declarations.} There is no conflict of interests.

\bye